\documentclass[12pt]{article}
\usepackage{amsfonts}
\usepackage{amsmath}
\usepackage{epsfig}

\title{On Area Growth in Sol}
\author{Richard Evan Schwartz \thanks{Supported by N.S.F. Grant DMS-1807320}}

\newtheorem{theorem}{Theorem}[section]

\newtheorem{lemma}[theorem]{Lemma}

\newtheorem{corollary}[theorem]{Corollary}

\def\startproof{{\bf {\medskip}{\noindent}Proof: }}

\def\endproof{$\spadesuit$  \newline}

\def\R{\mbox{\boldmath{$R$}}}% 

\begin{document}
\maketitle
\begin{abstract}
  Let Sol be the $3$-dimensional solvable Lie group
  whose underlying space is $\R^3$ and whose
  left-invariant Riemannian metric is given by
  $$e^{-2z} dx^2 + e^{2z} dy^2 + dz^2.$$
  We prove that the sphere of radius $r$
  in Sol has surface area at most $20 \pi e^r$ provided
  that $r$ is sufficiently large. This estimate is
  sharp up to a factor of $10$.
\end{abstract}

\section{Introduction}

Sol is one of the $8$ Thurston geometries [{\bf Th\/}],
the one which uniformizes torus bundles which fiber
over the circle with Anosov monodromy.   Sol has
been studied in various contexts:  coarse geometry
[{\bf EFW\/}],  [{\bf B\/}];  minimal surfaces
[{\bf LM\/}] (etc.);
its geodesics [{\bf G\/}], [{\bf T\/}], [{\bf K\/}],
[{\bf BS\/}]; connections to Hamiltonian systems
[{\bf A\/}], [{\bf BT\/}]; and finally virtual reality
[{\bf CMST\/}].  Our paper
 [{\bf CS\/}] has a more
extensive discussion of these many references.

In [{\bf CS\/}], Matei Coiculescu and I
give an exact characterization of
which geodesic segments in Sol are length minimizers,
thereby giving a precise description of the cut locus
of the identity in Sol.  As a consequence, we proved
that the metric spheres in Sol are topological spheres,
smooth away from $4$ singular arcs.  We will
summarize the characterization in \S 2 and
explain the main ideas in the proof in \S 4.

Ian Agol recently pointed out to me that our
exact characterization of the cut locus in
Sol might help us determine the growth rate
in balls in Sol.
Eryk Kopczy\'nski pointed out a rather
easy calculation that $V_r \leq C r^2 e^r$
for some constant $C$, and this implies that
Sol has volume entropy 1.
See \S \ref{eryk}.  Marc Troyanov recently
showed me a preprint with the estimate
$V_r<8 r e^r$.  We will get finer information.

To state our main result, we normalize the metric
in Sol so that it is:
\begin{equation}
  e^{-2z}dx^2 + e^{2z}dy^2+dz^2.
\end{equation}
In this metric
the planes $X=0$ and $Y=0$ have sectional curvature $-1$.
Let ${\cal S\/}_r$ denote any metric sphere of radius
$r$ in Sol.  Let $A_r$ denote the area of 
${\cal S\/}_r$ with respect to the Riemannian metric in Sol.

\begin{theorem}
  \label{main}
  $A_r<20 \pi e^r$ provided that $r$ is
  sufficiently large.
\end{theorem}

\noindent
    {\bf Remarks:\/} \newline
    (i) Given the $2$-to-$1$ locally-area-decreasing projection from
     ${\cal S\/}_r$ onto the hyperbolic disk of radius $r$, we have
     $A_r>4\pi(\cosh(r)-1) \approx 2\pi e^r$.  Thus, our estimate
is sharp up to a factor of $10$.
    \newline
(ii) Our analysis does not give
 an effective estimate on what
``sufficiently large'' means.
However, we notice that all
the relevant quantities seem to stabilize
pretty quickly: one sees all the phenonena
already when looking at a sphere of radius
$8$ in Sol. \newline

The basic idea of the proof is to bound the
projections of the sphere into the coordinate
planes. Let $\Pi_X$ denote the coordinate
plane $X=0$. Let $\eta_X$ denote the
projection of Sol into $\Pi_X$.
Let $A_{X,r}$ denote the area of
$\eta_X({\cal S\/}_r)$.  Let
$N_{X,r}$ denote the smallest integer such that
the map $\eta_X: {\cal S\/}_r \to \Pi_X$ is
at most $N_{X,r}$-to-$1$. We make all the
same definitions with $Y$ and $Z$ in place of $X$.
We also make a more refined definition for $Z$.
Let ${\cal S\/}_{r,k}$ denote the subset
where $\eta_Z$ is $k$-to-$1$,
We let $A_{Z,r,k}$ denote the area of
$\eta_Z({\cal S\/}_{r,k})$.
In \S \ref{projboundproof}
we prove the following result.

\begin{lemma}[Projection]
  \label{projbound}
For any $\epsilon>0$ and any $\theta \in (0,1)$ we may take $r$ sufficiently
large so that
$$A_r<\frac{(N_{X,r} A_{X,r}+N_{Y,r} A_{Y,r})}{\theta-\epsilon}+
\frac{1}{\sqrt{1-\theta^2}-\epsilon} \sum_{k=1}^{\infty} k A_{Z,k,r}. $$
\end{lemma}

\noindent
{\bf Remark:\/}
A similar result would be true for any surface
in Sol which is almost everywhere smooth, but the formula
in general would not be quite as good.  We use special
properties of ${\cal S\/}_r$ to get the formula above.
\newline

We then establish the following 
{\it projection estimates:\/}

\begin{enumerate}
  \item $A_{X,r}=A_{Y,r} =    2 \pi (\cosh(r)-1)<\pi e^r$.
  \item $N_{X,r}=N_{Y,r}=2$.
  \item $A_{Z,r}<16^*e^r$ for sufficiently large $r$.
  \item $N_{Z,r}=4$ for sufficiently large $r$.
  \item $A_{Z,k,r}<0^* e^r$ for $k=3,4$.
\end{enumerate}
The number $\zeta^*$ is a number we can make as close
as we like to $\zeta$ by taking $r$ sufficiently large.
Estimates 4 and 5 together combine to say that
the projection $\eta_Z$ is essentially $2$-to-$1$, because
the set where it is either $3$-to-$1$ or $4$-to-$1$ has
negligible area in comparison to $e^r$.
When we apply the result in Lemma \ref{projbound} we see that,
for $r$ large,
\begin{equation}
\label{bigsum}
A_r < \min_{\theta \in (0,1)} \bigg(\frac{4^* \pi}{\theta} + 
\frac{32^*}{\sqrt{1-\theta^2}}\bigg) e^r<20 \pi e^r.
\end{equation}
The minimizer is quite close to $\theta=3/5$ and
the minimum is about $60.93 e^r$.

Projection Estimates 1 and 2 are  straightforward
given our description of the Sol spheres.  Projection
Estimate 3 relies on an analysis of an ODE studied
in [{\bf CS\/}] and some easy asymptotic results about
elliptic functions.  After giving our upper bound we
will explain, a bit sketchily, why
$$A_{Z,r}>(2/1^*) e^r$$ once $r$ is
sufficiently large.  We do this to point out that our
upper bound is fairly tight.  When we plug in this
smaller estimate into the Projection formula above,
we get a bound of about $7\pi e^r$.  Ths represents a
kind of absolute limit to the strength of our method.

The hard
work in the paper involves dealing with
Projection Estimate 4, even though the result
is clear from the computer plots
such as Figure 5.4, involves a careful asymptotic
study of the ODE just mentioned.  Projection Estimate 5
comes out as a byproduct of our analysis.
\newline

This paper is organized as follows.
\begin{itemize}
\item In \S 2 we introduce some preliminary
  material and in particular recall the Main
  Theorem from [{\bf CS\/}].  We use this
  result to prove the Projection Lemma, though
  we also point out that one can prove the
  Projection Lemma knowing a much softer
  result about the Sol spheres.

\item In \S 3, we prove Projection Estimates 1 and 2.

\item In \S 4, we give details about the proof
  of the Main Theorem in [{\bf CS\/}].  These details
  are needed for the proof of Projection
  Estimates 3 -- 5.

\item In \S 5 we give more information about the central
  ordinary differential equations which arise in
  [{\bf CS\/}].

\item In \S 6 we prove Projection Estimates 3 -- 5
modulo the detail that a certain curve in the
plane is smooth and regular except at a single cusp.
See Figure 5.3.  Proving this result,
which we call the Embedding Theorem, turns out to
be a fight with the OEDs introduced in \S 5.

\item In \S 7 we prove the Embedding Theorem
  modulo a detail which we call the
  Monotonicity Lemma, a statement about
  the ODE from \S \ref{diffeq}.
  
\item in \S 8-9 we prove the Monotonicity
  Lemma.  This is where all the ODE calculations come in.
\end{itemize}

I would like to thank Ian Agol,  Matei Coiculescu,
Justin Holmer, Anton Izosimov, Boris Khesin,
Eryk Kopczy\'nski, Mark Levi, Benoit Pausader, Pierre Pansu,
and Marc Troyanov for helpful discussions
concerning this paper.
I would also like to
acknowledge the support of the Simons
Foundation, in the form of a 2020-21 Simons
Sabbatical Fellowship, and also the
support of the Institute for Advanced Study, in the
form of a 2020-21 membership funded by a grant
from the Ambrose Monell Foundation.

\newpage

\section{Preliminaries}

\subsection{Basic Properties of Sol}

The underlying space for Sol is $\R^3$.  The metric is:
\begin{equation}
  e^{-2z} dx^2 + e^{2z} dy^2 + dz^2.
  \end{equation}
The group law on Sol is
\begin{equation}
  \label{grouplaw}
  (x,y,z) * (a,b,c)= (e^z a+x,e^{-z} b+y,c+z).
\end{equation}
Left multiplication is an isometry.
We identify $\R^3$ with the Lie algebra of Sol in the obvious way.
(See [{\bf CS\/}, \S 2.1] if this does not seem obvious.)

Sol has $3$ interesting foliations.
\begin{itemize}
\item The XY foliation is by (non-geodesically-embedded) Euclidean planes.
 \item The XZ foliation is by geodesically embedded hyperbolic planes.
 \item The YZ foliation is by geodesically embedded hyperbolic planes.
\end{itemize}

The complement of the union of the
two planes $X=0$ and $Y=0$ is a union of
$4$ {\it sectors\/}.
One of the sectors, the
{\it positive sector\/}, consists of vectors
of the form $(x,y,z)$ with $x,y>0$. 
The sectors are permuted by the Klein-4
group generated by isometric reflections in the
planes $X=0$ and $Y=0$.  The Riemannian exponential
map $E$ preserves the sectors.
Usually, this symmetry will allow 
us to confine our attention to the positive sector.
\newline
\newline
{\bf Notation:\/}
For each $W \in \{X,Y,Z\}$, the plane $\Pi_W$ is given by $W=0$
and the map $\eta_W: {\rm Sol\/} \to \Pi_W$ is the projection
onto $\Pi_X$ obtained by just dropping the $W$ coordinate.
We also let $\pi_Z$ denote projection onto the $Z$-axis.
Thus, $\pi_Z(x,y,z)=z$.

\subsection{Properties of the Hyperbolic Slices}

We discuss our results for $\Pi_Y$.
There are analogous results for $\Pi_X$.
Here is a basic property of the hyperbolic slices
in Sol. The map $F(x,0,z)= (x,e^z)$
converts the metric in $\Pi_Y$ to the standard
hyperbolic metric in the upper half plane, namely
$$(dx^2+dy^2)/y^2.$$

\begin{lemma}
\label{fact1}
In $\Pi_Y$, the points $(0,0,0)$ and $(d_r,0,0)$
are connected by a geodesic segment of length $r$ when
$d_r=e^{r/2}-e^{-r/2}$.
\end{lemma}

\startproof
Here $d_r=2\sinh(r/2)$.
Let $F$ be the transformation from $\Pi_Y$ to the standard
upper half plane model.
We have
$F(0,0,0)=(0,1)$ and
$F(d_r,0,0)=(d_r,1)$.
As is well known, the distance between these points in the
standard hyperbolic metric is
$2\sinh^{-1}(d_r/2)=2(r/2)=r.$
\endproof

\begin{lemma}
\label{fact2}
In $\Pi_Y$, the point $(x,0,z)$ lies in the disk of
radius $r$ centered at $(0,0,0)$ only if
$|x| \leq (e^r-e^{-r})/2$.
\end{lemma}

\startproof
We use the transformation $F$ again.
Looking in the standard upper half plane model,
the disk we are interested in, $D_r$, is
centered at $(0,1)$ and has radius $r$.
The two points $(0,e^{-r})$ and $(0,e^r)$
lie in the boundary of $D_r$. Hence
$D_r$ has Euclidean radius 
$(e^r-e^{-r})/2$.
\endproof

\subsection{The Disk Lemma}
\label{disk0}

Here we recall a result from topology.
This result will be useful, in \S 5, when
we prove Projection Estimate 4.

\begin{lemma}[Disk]
  Let $\Delta \subset \R^2$ be a disk.
Let $h: \Delta \to \R^2$ be a map which is a local
diffeomorphism on the interior such that
$h(\partial \Delta)$ is a
  piecewise smooth curve having finitely
  many self-intersections.
Given $p \in \R^2-h(\Delta)$, the number of
preimages $h^{-1}(p)$ equals the unsigned number of
times $h(\partial \Delta)$ winds around $p$.
\end{lemma}

\startproof
This is a well-known result.
Here we sketch the
proof.  Without loss of generality, we can
assume that $\Delta$ is the unit disk in $\R^2$ and
$h(0,0) \not = p$.  Let $\Delta_s$ denote the
disk of radius $r$ centered at $(0,0)$.
Also, we can assume that $h$ is orientation
preserving in the interior of $\Delta$.
Let $f(s)$ denote the number of times
$h(\Delta_s)$ winds around $p$.   For $s$ near
$0$, we have $f(s)=0$.  The function
$f$ changes by $\pm 1$ each time
$\Delta_s$ crosses a point of $f^{-1}(p)$.
The sign is always the same because
$h$ is orientation preserving.
Hence the number of points
in $f^{-1}(p)$ equals $f(1)$, up to sign.
\endproof

\subsection{Elliptic Functions}

Many of the quantities associated to Sol are expressed in
terms of elliptic integrals.  Our functions ${\cal K\/}$ and
${\cal E\/}$ are precisely {\tt EllipticK\/} and {\tt EllipticE\/} in
Mathematica [{\bf W\/}].
\newline
\newline
{\bf Basic Definition:\/}
The complete elliptic functions of the first and second kind are given by
\begin{equation}
  \label{ellip}
  {\cal K\/}(m) = \int_0^{\pi/2} \frac{d\theta}{\sqrt{1-m\sin^2\theta}}, \hskip 20 pt
  {\cal E\/}(m)=\int_0^{\pi/2} \sqrt{1-m \sin^2(\theta)}\ d\theta.
\end{equation}
The first integral has domain $m \in [0,1)$ and the second has domain $m \in [0,1]$.
  \newline
  \newline
  {\bf Differential Equations:\/}
  These functions satisfy the following differential equations.
  For a proof see any textbook on elliptic functions.
  \begin{equation}
\label{diff}
\frac{d {\cal K\/}}{dm}=\frac{(m-1) {\cal K\/} + {\cal E\/}}{2m-2m^2}, \hskip 30 pt
\frac{d {\cal E\/}}{dm}=\frac{- {\cal K\/} + {\cal E\/}}{2m}.
  \end{equation}

  \noindent
  {\bf AGM Identity:\/}
  We have the following classic identity.
\begin{equation}
\label{agmid}
{\cal K\/}(m)=\frac{\pi/2}{{\rm AGM\/}(\sqrt{1-m},1)}, \hskip 30 pt
m \in (0,1).
\end{equation}
See [{\bf BB\/}] for a proof.
\newline
\newline
{\bf Asymptotics:\/}
It follows directly from the definition and an elementary integral that
\begin{equation}
  \label{ASY0}
 {\cal E\/}(1)=1.
 \end{equation}
We also have the following:
\begin{equation}
  \label{ASY}
  \bigg|{\cal K\/}(m) + \frac{1}{2} \log\bigg(\frac{1-m}{16}\bigg)\bigg|<
  \frac{1-m}{8} \times  \log\bigg(\frac{1-m}{16}\bigg).
\end{equation}
Both sides tend to $0$ as $m \to 1$.
This inequality comes from the second inequality (the upper bound) in
Inequality 19.9.2 of the Digital Library of Mathematical Functions:
$$1+\frac{(k')^2}{8} <\frac{{\cal K\/}(k)}{\log(4/k')} < 1+\frac{(k')^2}{4}.$$
Here $m=k^2$ and $1-m=(k')^2$.  Making the substitution of $m$ and $1-m$
for $k$ and $k'$ and then rearranging, we get
Equation \ref{ASY}.

\subsection{The Hamiltonian Flow}

Let $G={\rm Sol\/}$.
Let $S_1 \subset \R^3$ denote the unit sphere.
at the origin in $G$.  Given a unit speed geodesic
$\gamma$, the tangent vector $\gamma'(t)$ is part of a
left invariant vector field on $G$, and we
let $\gamma^*(t) \in S_1$ be the restriction of
this vector field to $(0,0,0)$.  In terms of
left multiplication on $G$, we have the formula
\begin{equation}
\label{star}
\gamma^*(t)=dL_{\gamma(t)^{-1}}(\gamma'(t)).
\end{equation}

It turns out that
$\gamma^*$ satisfies the following
differential equation.
\begin{equation}
\label{GraysonStructureField}
  \frac{d\gamma^*(t)}{dt}=\Sigma(\gamma^*(t)),
  \hskip 30 pt \Sigma(x,y,z)=(+xz,-yz,-x^2+y^2).
\end{equation}
This is explained one way in [{\bf G\/}] and 
another way in  [{\bf CS\/}, \S 5.1].
(Our formula has a different sign than Grayson's, because our group law
correspondingly differs by a sign.)
This system in Equation \ref{GraysonStructureField} is really
just geodesic flow on the
unit tangent bundle of Sol, viewed in a left-invariant reference frame.

Let $F(x,y,z)=xy$.  The flow lines
of $\Sigma$ lie in the level sets of $F$,
and indeed $\Sigma$ is the Hamiltonian flow
generated by $F$.
Most of the level sets of $F$ are closed loops.
We call these {\it loop level sets\/}.
With the exception of the points in the planes $X=0$ and $Y=0$,
and the points $(x,y,0)$ with $|x|=|y|=\sqrt 2$, the
remaining points lie in loop level sets.

Each loop level set $\Theta$ has an associated {\it period\/} $L=L_{\Theta}$,
which is the time it takes a flowline -- i.e., an integral curve -- in $\Theta$
to flow exactly once around.  Equation \ref{periodXX} below
gives a formula. We can compare $L$ to the
length $T$ of a geodesic segment $\gamma$ associated to a flowline
that starts at some point of $\Theta$ and flows for time $T$.
We call $\gamma$
{\it small\/}, {\it perfect\/}, or {\it large\/} according
as $T<L$, or $T=L$, or $T>L$.
In [{\bf CS\/}, \S 5] we prove the following result:
\begin{theorem}
  Suppose $(x,y,z) \in S^2$ lies in a loop level set.
  Let $\alpha=\sqrt{|xy|}$.  Then the period of the loop level
  set containing $(x,y,z)$ is
  \begin{equation}
\label{periodXX}
L_\alpha=\frac{\pi}{{\rm AGM\/}(\alpha,\frac{1}{2}\sqrt{1+2\alpha^2})}=
\frac{4}{\sqrt{1+2\alpha^2}} \times {\cal K\/}\bigg(\frac{1-2\alpha^2}{1+2\alpha^2}\bigg)
\end{equation}
\end{theorem}
The second identity is Equation \ref{agmid}.
Using Equation \ref{ASY} we see that the difference between
$L_{\alpha}$ and $-4\log(\alpha/2)$ tends to $0$ as $\alpha \to 0$.

Each vector $V=(x,y,z)$ simultaneously corresponds to two
objects:
\begin{itemize}
\item The flowline $\phi_V$ which starts at $V/\|V\|$ and
  goes for time $\|V\|$.
\item The geodesic segment $\gamma_V=\{E(tV)|\ t \in [0,1]\}.$
\end{itemize}
Given a vector $V=(x,y,z)$ we define
\begin{equation}
  \mu(V)={\rm AGM\/}(\sqrt{xy},\frac{1}{2}\sqrt{(|x|+|y|)^2+z^2}).
\end{equation}
We call $V$ {\it small\/}, {\it perfect\/}, or {\it large\/} according
as $\mu(V)$ is less than, equal to, or greater than $\pi$.  In view
of Equation \ref{periodXX} here is what this means:
\begin{itemize}
\item If $\gamma_V$ lies in the plane $X=0$ or $Y=0$ then $V$ is small
  because $\mu(V)=0$.
\item If $V=(x,y,0)$ where $|x|=|y|$ then $V$ is small, perfect, or
  large according as $|x|<\pi$, $|x|=\pi$ or $|x| \geq \pi$.
\item In all other cases, $V/\|V\|$ lies in a loop level
  set, and $V$ is small, perfect, or large according as
  $\mu(V)<\pi$, $\mu(V)=\pi$, or $\mu(V)>\pi$.
\end{itemize}

\subsection{The Main Result}

Now we recall the main result from [{\bf CS\/}].

\begin{theorem}
  \label{main2}
  Given any vector $V$, the geodesic segment $\gamma_V$ is a distance
  minimizing geodesic if and only if $\mu(V) \leq \pi$.  That is,
  $\gamma_V$ is distance minimizing if and only if $V$ is small or perfect.
  Moreover, if $V$ and $W$ are perfect vectors then
  $E(V)=E(W)$ if and only if $V=(x,y,z)$ and $W=(x,y,\pm z)$.
  \end{theorem}

Theorem \ref{main} identifies the cut locus of the identity in
Sol with the set of perfect vectors.  The Riemannian
exponential map $E$ is a global diffeomorphism on the
set of small vectors.  Also, $E$ 
is generically $2$-to-$1$ on the set of perfect vectors.
We will explain this last fact below.

Theorem \ref{main} leads to a good  description of 
the Sol  metric sphere ${\cal S\/}_r$ of radius $r$.
Let $S_r$ denote the Euclidean sphere of radius $r$
centered at the origin of $\R^3$.  Let
\begin{equation}
  S'_r=\mu^{-1}[0,\pi] \cap S_r.
\end{equation}
The space $S_r'$ is a $4$-holed sphere.  The
boundary $\partial S_r'$, a union of
$4$ loops, is precisely the set of
perfect vectors contained in $S_r$.
Each of these loops is homothetic to one
of the loop level sets on the unit sphere $S^2$.
The Klein-4 symmetry explains why there are
$4$ such loops.

It follows from the Main Theorem that
${\cal S\/}_r=E(S'_r)$ and that
$E$ is a diffeomorphism when restricted to $S'_r-\partial S'_r$.
On $\partial S'_r$, the map $E$ is a
$2$-to-$1$ folding map which
identifies partner points within each component.  Thus, we see that
${\cal S\/}_r$ is obtained from a $4$-holed sphere by gluing
together each boundary component (to itself) in a $2$-to-$1$ fashion.
This reveals
${\cal S\/}_r$ to be a topological sphere which is smooth away from
the set $E(\partial S'_r)$.  We also prove that 
the singular set $E(\partial S_r')$ consists of $4$ arcs of hyperbolas, all
contained in $\Pi_Z$.
\newline
\newline
{\bf The Lunar Principle:\/}
Given a unit normal vector $V$ to ${\cal S\/}_r$ at a 
smooth point, we let $V_*$ denote the left translate
of $V$ to the origin.  Let $N_r$ denote the set of all
such vectors $V_r$.  Given the nature of the loop level
sets, we have the following corollary of Theorem \ref{main2}.
For any $\epsilon>0$ there is some $R$ such that
$N_r$ is contained in the $\epsilon$-tubular neighborhood
of $\Pi_X \cup \Pi_Y$ provided that $r>R$.  We call this
the {\it Lunar Principle\/} because $\Pi_X \cup \Pi_Y$
intersects the unit sphere in a union of $4$ spherical
lunes.  One does not really need the full force of
Theorem \ref{main2} to deduce the Lunar Principle:
A long geodesic tangent to
a unit vector that is far from $\Pi_X \cup \Pi_Y$
makes a corkscrew-like pattern and is quite far
from distance minimizing.

\subsection{Proof of the Projection Lemma}
\label{projboundproof}

We call a map $\eta$ between surfaces
$\theta$-{\it good\/} if 
$$\frac{{\rm area\/}(\eta(S))}{{\rm area\/}(S)} \geq \delta$$
for any measurable subset $S$ in the domain.
Mostly we are interested in the case when the domain
is a smooth surface in Sol and the range is one of the
coordinate planes in Sol.  However, in the first result,
we will consider planar surfaces in $\R^3$.  The same
projections $\eta_X,\eta_Y,\eta_Z$ make sense as
projections in $\R^3$.

\begin{lemma}
\label{proj0}
Let $\Theta_X,\Theta_Y,\Theta_Z$ be 
positive numbers with
$\Theta_X^2+\Theta_Y^2+\Theta_Z^2=1$.  Let $\Pi$
be any plane in $\R^3$.  Then there is
an $I \in \{X,Y,Z\}$ such that
the $\eta_I$ is $\Theta_I$-good.
\end{lemma}

\startproof
It follows from the familiar fact that $\|V \times W\|$
computes the area of the parallelogram spanned by two
vectors $V,W \subset \Pi$, and from the Pythagorean Theorem,
that there are $3$ non-negative numbers
$r_X,r_Y,r_Z \geq 0$ so that
$r_X^2+r_Y^2+r_Z^2=1$, and 
\begin{equation}
\label{pyth}
A_{\Pi}(S)=r_X A_X(S)+r_Y A_Y(S)+ r_Z A_Z(S).
\end{equation}
Here $S \subset \Pi$ is any measurable set and
$A_X(S)$ is the area of $\Pi_X(S)$, etc.
If our claim is false then $\Theta_I<r_I$ for all
$I \in \{X,Y,Z\}$.  But then
$$1=\Theta_X^2+\Theta_Y^2+\Theta_Z^2<r_X^2+r_Y^2+r_Z^2=1,$$
and we have a contradiction.
\endproof

Now we move the discussion to Sol.

\begin{lemma}
\label{proj1}
Let $\Theta_X,\Theta_Y, \Theta_Z$ be 
positive numbers with
$\Theta_X^2+\Theta_Y^2+\Theta_Z^2=1$.
Let $\Sigma$ be a smooth surface in Sol.
Let $p \in \Sigma$ be any point.
Then for any $\epsilon>0$ there is a
sufficiently small neighborhood $U$ about
$p$ and some index $I \in \{X,Y,Z\}$ such that
$\eta_I$ is $\Theta_I$-good on $U$.
\end{lemma}

\startproof
Given that Sol is homogeneous, and that the
projections between parallel planes within the
same coordinate foliation are area preserving,
it suffices to prove our result when $p$ is
the origin in Sol.  But, in this case, the
metric on Sol agrees with the Euclidean metric
up to any given $\epsilon$ we like.  So, this
special case follows from Lemma \ref{proj0}
and the differentiability of $\Sigma$.
\endproof

Now let us apply the Lunar Principle to the sphere
${\cal S\/}_r$.  Given a smooth point
$p \in {\cal S\/}_r$, the corresponding
vector $N_{p,0}$ lies quite near
$\Pi_X \cup \Pi_Y$.  The reason is that
$N_{p,0}$ either lies in $\Pi_X \cup \Pi_Y$
or else in a loop level set of period
greater than $r$, and such loop level
sets lie near $\Pi_X \cup \Pi_Y$.
Therefore, given any $\epsilon>0$ we can
take $r$ large enough so that there is a
partition of the smooth points of ${\cal S\/}_r$
into $3$ measurable (or indeed piecewise smooth) regions 
$${\cal S\/}_r(I), \hskip 30 pt 
I \in \{X,Y,Z\}$$ with the following
properties:
\begin{itemize}
\item The projection $\eta_X: {\cal S\/}_r(X) \to \Pi_X$
is $(\theta-\epsilon)$ good.

\item The projection $\eta_Y: {\cal S\/}_r(Y) \to \Pi_Y$
is $(\theta-\epsilon)$ good.

\item The projection $\eta_Z: {\cal S\/}_r(Z) \to \Pi_Z$
is $(\sqrt{1-\theta^2}-\epsilon)$-good.
\end{itemize}
Since the non-smooth subset of ${\cal S\/}_r$ has
area $0$, the formula in the 
Projection Lemma follows immediately.

\subsection{A Weaker Bound on Volume}
\label{eryk}

This section is independent from the
rest of the paper. Here we present,
with minor modifications,
Eryk Kopczy\'nski's derivation of
a weaker volume growth bound that is
still sufficient to establish that
Sol has volume entropy 1.  I did not
try for optimal constants.

\begin{lemma}
  Suppose $\gamma$ is a geodesic in Sol
  having length $r$.  Let $r_1$ be the
  length of $\gamma$ that lies above the
  plane $Z=0$ and let $r_2$ be the length
  of $\gamma$ that lies below or in the plane $Z=0$.
  Then the endpoint $(x,y,z)$ of $\gamma$
  satisfies the bound 
  $|x| \leq e^{r_1}+r_2$ and $|y|<e^{r_2}+r_1$.
\end{lemma}

\startproof
We first consider two special cases.
If $\gamma$ stays above the plane $\Pi_Z$ then
the endpoint $(x,y,z)$ satisfies the bounds
$|x| \leq e^r$ and $|y| \leq r$.
Likewise, if $\gamma$ stays below the plane $\Pi_Z$ then
the endpoint $(x,y,z)$ satisfies the bounds
$|y| \leq e^r$ and $|x| \leq r$.
In general, one can break $\gamma$ into intervals
$\gamma_1,...,\gamma_k$, for some $k$, such that
each $\gamma_j$ satisfies one of the two
special cases just considered.  Adding up the
bounds from the special cases, we get
the result advertised in the lemma.
\endproof

Set $u=e^r$.  Note that $r<u$.
Every point $(x,y,z)$ in the ball of radius $r$ satisfies
$$|x|, |y| \leq u+r, \hskip 30 pt
|z| \leq r, \hskip 30 pt (|x|-r)(|y|-r) \leq u.$$
Let $\Omega_r$ be the set of points satisfying these inequalities.
For convenience we take $r \geq 1$
The volume of the part of $\Omega_r$ where $|x| \leq r+1$ is
bounded by
$$8r \times (r+1) \times (u+r) < 32 r^2 u.$$
Likewise, the  volume of the part of $\Omega_r$ where $|y| \leq r+1$ is
bounded by $32r^2 u$.
The volume of the part of the $\Omega_r$ where $|x|>r+1$ and $|y|>r+1$ is
$$8r \int_1^u \frac{u}{x}\ dx \leq 8ru \log(u)=8r^2 u.$$
Therefore $\Omega_r$ has volume at most $72r^2e^r$.
But the ball of radius $r$ is contained in $\Omega_r$.

\newpage

\section{The Hyperbolic Projections}

\subsection{A Picture}

Recall that $\eta_X: {\rm Sol\/} \to \Pi_X$ is the orthogonal
projection onto the plane $X=0$.
Figure 3.1 shows the projection of (part of) the positive sector
of the sphere ${\cal S\/}_5$ into the plane $\Pi_X$.
The smooth part of this sphere has a foliation by the
images of the loop level sets under the Riemannian exponential map $E$.
The grey curves are the projections of this foliation into $\Pi_X$.
The black line segment is the projection of the set of
singular points.

\begin{center}
\resizebox{!}{2.5in}{\includegraphics{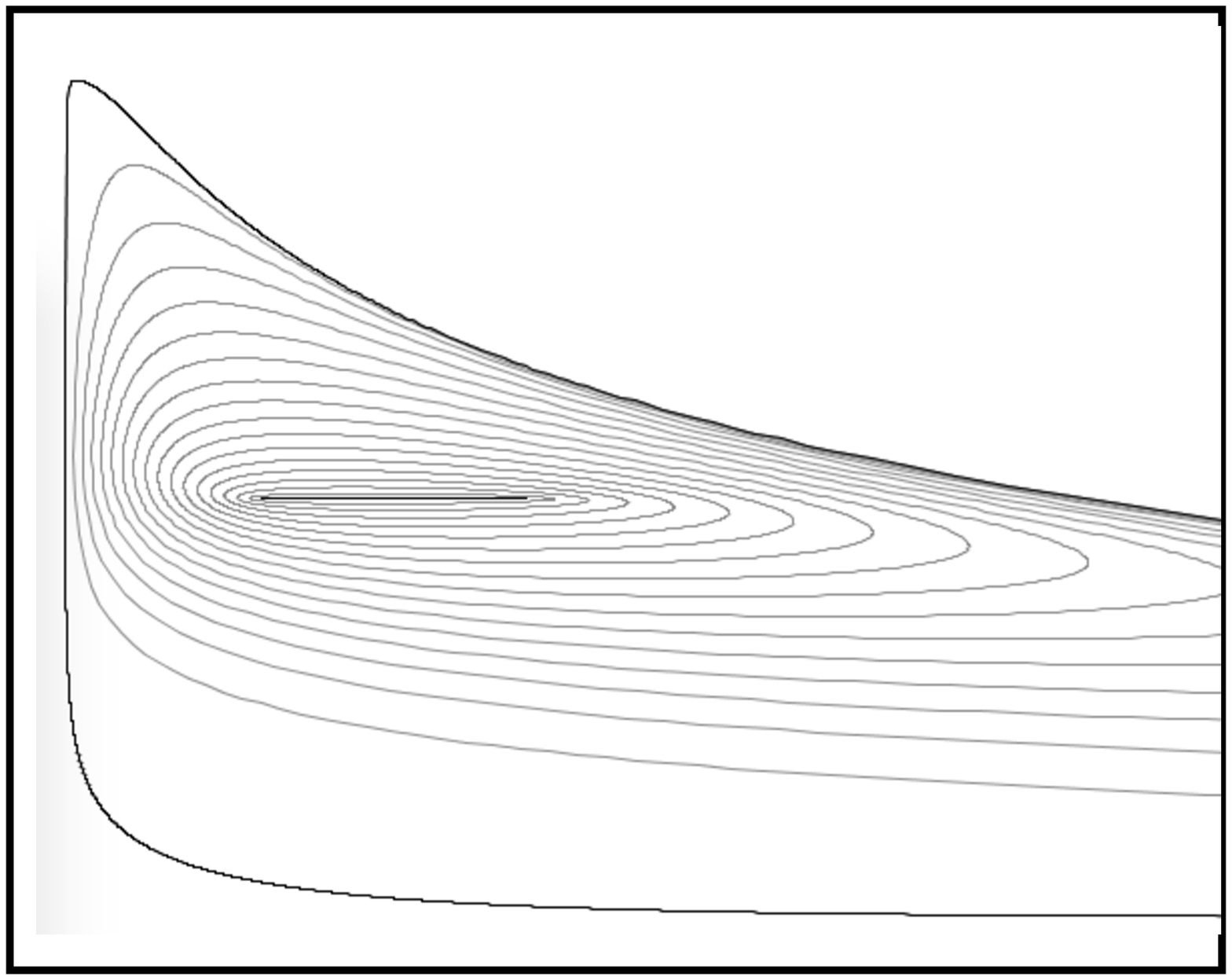}}
\newline
    {\bf Figure 3.1:\/} Projection into the plane $\Pi_X$.
\end{center}

 It appears
from the picture that the restriction of
$\eta_X$ to this sector is a homeomorphism onto its image.
We will prove this result below.
 For convenience we take $r>\pi \sqrt 2$.

\subsection{Area Bound}

In this section we prove Projection Estimate 1.
By symmetry, it suffices to prove the result
for the projection $\eta_X$ into the plane $\Pi_X$.
Let ${\cal S\/}_r^+$ denote the subset
of ${\cal S\/}_r$ consisting of points $(x,y,z)$ with $x \geq 0$.
Let ${\cal H\/}_r$ denote the hyperbolic disk of
radius $r$ contained in the plane $\Pi_X$ and centered at the origin.
All the points in the interior of ${\cal S\/}_r^+$
lie in the open positive sector.  Because $\Pi_X$ is
a totally geodesic plane in Sol, we have
$\partial {\cal S\/}_r^+= \partial {\cal H\/}_r$.
The following result immediately implies that
$A_{X,t}<\pi e^r$.

\begin{lemma}
  \label{hyp1}
  $\eta_X$ maps the interior of ${\cal S\/}_r^+$ into the interior
  of ${\cal H\/}_r$.
\end{lemma}

\startproof
If this is false, then there is a geodesic segment $\gamma$ of length $r$,
connecting $(0,0,0)$ to some point $p \in {\cal S\/}_r^+$ which
remains entirely in the positive sector except for its initial
point, $(0,0,0)$.  The projection map $\eta_X$ is distance
non-increasing, and locally distance decreasing on any curve
whose tangent vector is not in a plane of the form $X={\rm const.\/}$.
This means that $\eta_X(\gamma)$ is shorter than $\gamma$.
But then $\eta_X(\gamma)$ cannot reach the point
$\eta_X(p) \in  \partial {\cal H\/}_r$.
\endproof

\subsection{Multiplicity Bound}

In this section we prove Projection Estimate 2.
As above, it suffices to prove this result
for the projection $\eta_X$ into the plane $\Pi_X$.
To prove that $N_{X,r}=2$ it suffices, by symmetry,
to show that $\eta_X$ is an injective map from
${\cal S\/}_r^+$ to $\Pi_X$.  
The basic strategy is to show that
$\eta_X$ is locally injective.  We also know that
$\eta_X$ is the identity on the boundary of
${\cal S\/}_r^+$, which already lies in the plane $X=0$.
(It is the boundary of the hyperbolic disk on $\Pi_X$ of
radius $r$ centered at the origin.)  Our injectivity result
then follows from the Disk Lemma in \S 2.

For convenience we take $r>\pi\sqrt 2$ in the next result,
so that we don't have to discuss several cases. (The
sphere ${\cal S\/}_r$ is smooth for $r<\pi \sqrt 2$ and
has $4$ singular arcs for $r>\pi \sqrt 2$.)  The set
of smooth points of ${\cal S\/}_r$ is a union
of $4$ open ``punctured'' disks.  In each case, we
are removing an analytic arc from an open topological
disk and what remains is smooth.  Figure
4.1 shows (a portion of) the $\eta_X$-projections of the smooth points
of ${\cal S\/}_r$.

\begin{lemma}
  \label{hyp2}
  The differential $d\eta_X$ is injective at all the smooth points
  in the interior of ${\cal S\/}_r^+$
\end{lemma}

\startproof
Let $S_r$ denote the subset of the sphere of radius
$r$ centered at the origin in the Lie algebra.
As in the previous chapter, let
$S_r'$ denote the subset of $S_r$ consisting of
vectors which are either small or perfect.
Let $p \in S_r'$ be some point.  We think of $p$ as a
vector, so that $E(p) \in {\cal S\/}_r^+$.
Let $T_p$ be the
tangent plane to $S'_r$ at $p$.  Let $N_p$ be the
unit normal to $T_p$.    Since the perfect geodesic segments
are minimizers, the small geodesic segments are
unique minimizers without conjugate points. So, at the
corresponding points of $S'_r$, the differential
$dE_p$ is an isomorphism.  We just have to
show that $(1,0,0) \not \in dE_p(T_p)$.  We will
suppose that $(1,0,0) \in dE_p(T_p)$ and
derive a contradiction.

If $(1,0,0) \in dE_p(T_p)$, then
the first component of $dE_p(N_p)$ is $0$,
because $dE_p(N_p)$ and $dE_p(T_p)$ are
perpendicular.
Let $\gamma_p$ be the geodesic segment corresponding to $p$.
The vector $dE_p(N_p)$ is the unit vector tangent to
$\gamma_p$ at its far endpoint -- i.e., the endpoint
not at the origin.  This vector lies in the same
left invariant vector field as the endpoint $U_p$
of the flowline corresponding to $p$.  If the first
coordinate of $U_p$ is $0$, then the
entire flowline lies in the plane $X=0$. But then
$E(p) \in \partial {\cal S\/}_r^+$.
This is a contradiction.
\endproof

\begin{lemma}
\label{diff1}
  The map $\eta_X$ is locally injective at each
  singular point of ${\cal S\/}_r^+$.
\end{lemma}

\startproof
As we showed in [{\bf CS\/}], the singular
set in ${\cal S\/}_r$ consists of $4$ arcs
of hyperbolas, each contained in the plane
$\Pi_Z$.  Each of these arcs lies in the
interior of a different sector and is an
arc of a hyperbola.  These hyperbolas are
all graphs of functions.  The restriction
of $\eta_X$ to each hyperbola is therefore
injective.  We still need to see, however,
that $\eta_X$ is injective in neighborhoods
of these singular sets, and not just on
the singular sets.  There are two cases.
\newline
\newline
{\bf Case 1:\/}
Consider a point $p$ in the interior of
the singular set in ${\cal S\/}_r^+$.
By symmetry it suffices to consider
the case when $p$ is in the positive
sector.  The point $p$ lies in the plane
$\Pi_Z$ and has its first two coordinates positive.
There are exactly $2$ points $p_+,p_- \in S_r'$
such that $E(p_+)=E(p_-)=p$.  These points
have the form $p_1=(x,y,z)$ and $p_-=(x,y,-z)$.
We called such points {\it partners\/}.
In [{\bf CS\/}, Lemma 2.8] we showed that
$dE_{p_{\pm}}$ is non-singular.  This crucially
uses the fact that $p_{\pm}$ is a perfect
vector whose third coordinate is nonzero.
The same argument as in the previous lemma now
shows that the linear map $\eta_X \circ dE_{p_{\pm}}$
is an isomorphism from the tangent plane
$T_{p_+}$ to $\R^2$.  But then $\eta_X \circ E$ is a
diffeomorphism when restricted to
an open neighborhood $U_{\pm}$ of
$p_{\pm}$ in $S_r'$.
\begin{center}
\resizebox{!}{1.3in}{\includegraphics{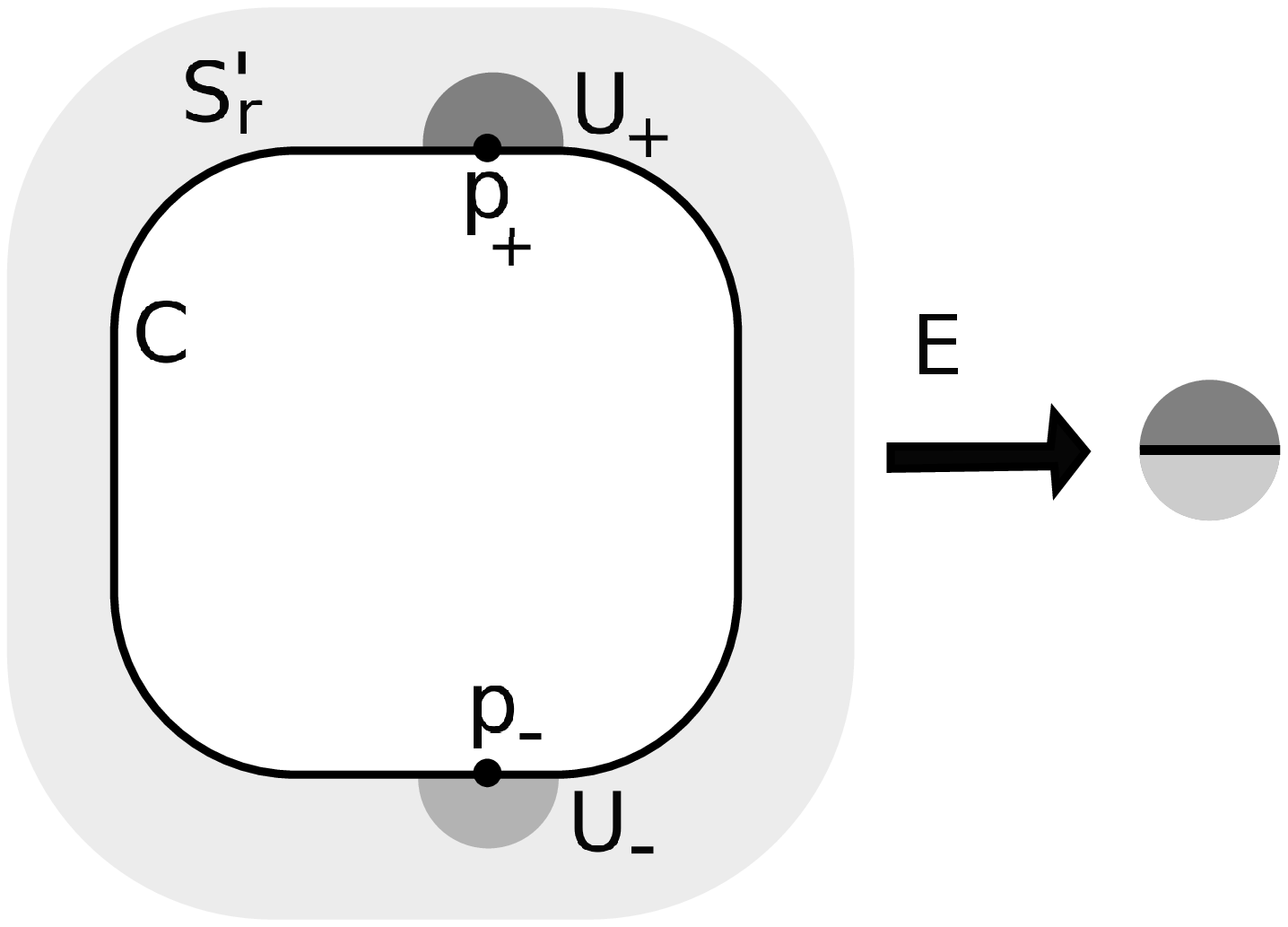}}
\newline
    {\bf Figure 3.2:\/} The neighborhoods $U_+$ and $U_-$.
\end{center}
The sets
$U_{\pm}$ are disks with some of
their boundary included.  The
portion of the included boundary consists of
the perfect vectors in $\partial S_r'$ near
$p_{\pm}$.  See Figure 3.2.

Let $C$ be the component of $\partial S_r'$ which
contains $p_+$ and $p_-$.  The image
$E(U_+-C)$ lies entirely below the plane
$\Pi_Z$ because the flowlines corresponding
to vectors in $U_+-C$ nearly wind the entirely
around their loop level set but omit a small
arc near $p_+$.   Likewise, the image
$E(U_--C)$ lies entirely above the
plane $\Pi_Z$.  Hence $\eta_X \circ E(U_+-C)$ and
$\eta_x \circ E(U_--C)$ are disjoint.
Combining this what we know, we see that
$\eta_X$ is a homeomorphism in a
neighborhood of $p \in {\cal S\/}_r$.
\newline
\newline
{\bf Case 2:\/}
Suppose that $p$ is one of the endpoints
of the singular set.  This case is
rather tricky to check directly.
Suppose that there is some other point
$q \in {\cal S\/}_n$ such that
$\eta_X(p)=\eta_X(q)$.  By Lemma \ref{hyp1}, the point $\eta_X(p)$ is
disjoint from the the hyperbolic
circle ${\cal S\/}_r^+ \cap \Pi_X$.
Hence $q$ lies in the interior of ${\cal S\/}_r^+$.
Since $\eta_X$ is injective on the singular set,
$q$ must be a smooth point.

Since $\eta_X(q)=\eta_X(p)$ and $p \in \Pi_Z$, we have
$q \in \Pi_Z$.  This means that $q$ corresponds to
some small symmetric flowline. The point $q$ is contained
in a maximal connected arc ${\cal A\/} \subset {\cal S\/}_r^+$
consisting entirely of points corresponding to
small symmetric flowlines.  One endpoint of
$\cal A$ is $p$. The other endpoint lies in 
the plane $X=0$.  The point $q$ lies somewhere
in the interior of $\cal A$.  The map
$\eta_X$ sends $\cal A$ into the line $X=Z=0$
and from Lemma \ref{diff1}, the restriction
of $\eta_X$ to $\cal A$ is locally injective.
But a locally injective map from an arc into a
line is injective.  This contradicts
the fact that $\eta_X(p)=\eta_X(q)$.
\endproof

Let $D$ denote the quotient $S_r^+/\sim$
where the equivalence relation $\sim$ glues together
partner points on the set of perfect vectors
in $S_r^+$.  The space $D$ is a topological
disk, and $h=\eta_X \circ E$ gives a map
from $D$ to $\Pi_X$. Combining
Lemmas \ref{hyp2} and \ref{diff1} we see that
$h$ is locally injective at each interior
point of $D$.  Moreover, $h(\partial D)$ is an
embedded loop, just the boundary of a hyperbolic
disk in $\Pi_X$.  By the Disk Lemma,
$h: D \to \eta_X$ is injective.
But $E$ is a bijection from
$D$ to ${\cal S\/}_r^+$.  Hence
$\eta_X: {\cal S\/}_r^+ \to \Pi_X$ is
injective, as desired.

This completes the proof that
$N_{X,r}=2$.

\newpage

\section{Details about the Cut Locus Theorem}
\label{solproof}

\subsection{Concatenation}
\label{concat}

In the next several sections, we outline the proof
of Theorem \ref{main2}.
Our exposition here is an abbreviated version
of what appears in [{\bf CS\/}].

Given a (finite) flowline $g$ we write
$g=a|b$ if $g$ is the concatenation of
flowlines $a$ and $b$.  That is,
$a$ is the initial part of $g$ and $b$ is the final part.
We call $g$ {\it symmetric\/} if the endpoints
of $g$ have the form $(x,y,z)$ and $(x,y,-z)$.

Let $\Lambda_g$ denote the endpoint of the geodesic
segment associated to $g$, when this geodesic segment
starts at the origin.  It follows from
left-invariance of the metric that
\begin{equation}
\Lambda_g = \Lambda_a * \Lambda_b.
\end{equation}
Since the third coordinates of elements of
Sol commute, we have
\begin{equation}
  \label{additive}
\pi_Z(\Lambda_g) = \pi_Z(\Lambda_a)+\pi_Z(\Lambda_b).
\end{equation}
Here $\pi_Z(x,y,z)=z$. More formally,
$\pi_Z$ is the quotient map from Sol to the quotient
${\rm Sol\/}/\Pi_Z$.  Here $\Pi_Z$ is not just
a Euclidean plane in Sol but also a maximal
normal subroup.
The integral form of Equation \ref{additive} is
\begin{equation}
\label{integralZ}
\pi_Z(\Lambda_g)=\int_0^T z(t)\ dt.
\end{equation}
Here we have set $g=(x,y,z)$, and
$T$ is the total time that $g$ takes to
get from start to finish.

These equations have a variety of consequences,
which we work out in detail in [{\bf CS\/}, \S 2].
\begin{enumerate}
\item If $g$ is a small flowline then $g$ is symmetric if and only if $\pi_Z(\Lambda_g)=0$.
Moreover, the geodesic segment corresponding to a small symmetric flowline only
intersects $\Pi_Z$ at its endpoints.
\item If $g$ is a perfect flowline then $\pi_Z(\Lambda_g)=0$.  This follows from
  the fact that $g=a|b$ where $a$ and $b$ are both small symmetric.
\item If $V_{\pm}=(x,y,\pm z)$, then $V_+$ is perfect if and only if $V_-$ is perfect.
Furthermore $E(V_+)=E(V_-)$.  This is because the corresponding flowlines
$g_+$ and $g_-$ can be written as $g_+=a|b$ and $g_-=b|a$ where
$a$ and $b$ are both small symmetric.  But then $\Lambda_a$ and $\Lambda_b$
are horizontal translations in Sol and hence commute.
Hence $\Lambda_{g_+}=\Lambda_{g_-}$.  We call
$V_+$ and $V_-$ {\it partners\/}.
\item Suppose $V_1$ and $V_2$ are perfect vectors such that
$V_1/\|V_1\|$ and $V_2/\|V_2\|$ lie in the same loop level set.
Let $E(V_i)=(a_i,b_i,0)$.
We call $\sqrt{a_ib_i}$ the {\it holonomy\/} of $V_i$.
Letting $g_1$ and $g_2$ be the corresponding flowlines,
we can write $g_1=a|b$ and $g_2=b|a$ where $a$ and $b$ are
both small. But then $\Lambda_{g_1}=(a_1,b_1,0)$ and
$\Lambda_{g_2}=(a_2,b_2,0)$ are conjugate in Sol.
This gives $a_1b_1=a_2b_2$. Hence $V_1$ and $V_2$ have
the same holonomy. 
\item Given $V=(x,y,z)$ we define
$\sigma(V)=y/x$.  We prove that if $V$ is a perfect
vector, then $\sigma(E(V))=1/\sigma(V)$.  We call this
the Reciprocity Lemma.  The proof is a more subtle working
out of the consequences of the conjugacy idea discussed
in Item 4.
\end{enumerate}

\subsection{Outline of the Proof}
\label{proofoutline}

With these preliminaries out of the way, we turn
directly to the proof of Theorem \ref{main2}.
Item 3 in \S \ref{concat}
shows that the perfect geodesic segments corresponding to
vectors of the form $(x,y,z)$ where $z \not =0$ are not
unique distance minimizers.  It also follows from
Item 3 that perfect geodesics segments corresponding
to vectors of the form $(x,y,0)$ have conjugate points.
Hence, large geodesic segments cannot be distance
minimizers.  This essentially proves half of
Theorem \ref{main2}.

The second half of Theorem \ref{main2}, the converse, says that
a small or perfect geodesic segment is a distance minimizer.
Since every small geodesic segment is
contained in a perfect geodesic segment, it suffices to
prove that perfect geodesic segments are
distance minimizers.  

We first prove [{\bf CS\/}, Corollary 2.10]:
The map $E$ is injective on the set of
perfect vectors with positive coordinates.
This step has $2$ ideas.
We first show (following [{\bf G\/}]) that the
holonomy is a monotone function of the loop
level set.  So, if $E(V_1)=E(V_2)$ then 
$V_1/\|V_1\|$ and $V_2/\|V_2\|$ lie in the
same loop level set.  We also
have $\sigma(V_1)=\sigma(V_2)$, by
Item 5 above. This forces $V_1=V_2$.

We finish the proof by showing that
if $V$ is perfect and $W$ is
small then it is impossible for $E(V)=E(W)$.
This is really the heart of [{\bf CS\/}].
The argument involves the system of nonlinear
ODEs we introduce in \S \ref{diffeq}.
It will turn out that the argument in this
paper involves a deeper study of these
same ODEs.

Let us go back to the argument.
By symmetry, we can restrict ourselves to the
case when $V$ and $W$ both lie in the positive
sector.  Let $M$ and $\partial M$ respectively denote the
set of small and perfect vectors.
We show that $E(\partial M)$
is contained in a subset $\partial N \subset \Pi_Z$.
The boundary of $\partial_N$, which we denote
by $\partial_0 N$, is the graph of a smooth
function in polar coordinates.  The yellow
region in Figure 3.1 shows part of the portion of
$\partial N$ that lies in the positive sector.
See Figure 5.1 for an expanded view.
There $3$ symmetrically placed components in
the other sectors which we are not showing.

\begin{center}
\resizebox{!}{2.5in}{\includegraphics{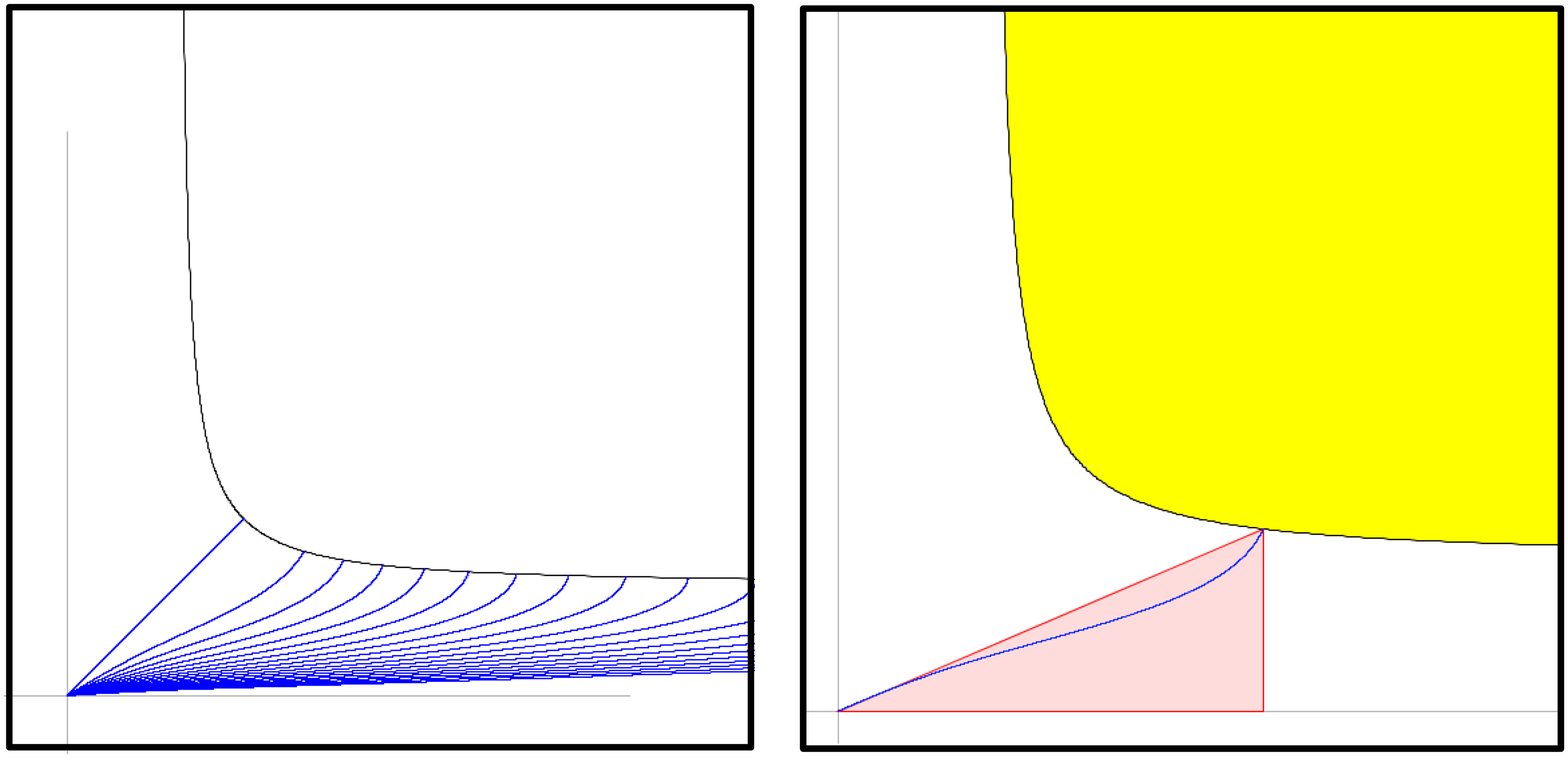}}
\newline
    {\bf Figure 3.1:\/} $\partial_0N_+$ (black), $\partial N_+$ (yellow),
    $\Lambda_L$ (blue), and $\Delta_L$ (red).
\end{center}

If we suppose that $W$ is small and $E(W)=E(V)$ then the
flowline corresponding to $W$ must be small symmetric.
We can arrange all the small symmetric flowlines in a
given loop level set into two curves.  One of the
curves corresponds to small symmetric flowlines whose
initial endpoint has positive $Z$-coordinate.
Given the loop level set of period $L$ in the positive sector,
we let $\Lambda_L$ denote the image, under $E$,
of the corresponding vectors.   The blue curves
in Figure 3.1 show $\Lambda_L$ for various
choices of $L$.

On the right side of Figure 3.1 we focus on $\Lambda_5$.  
We also draw the right triangle $\Delta_5$ whose
endpoints are the endpoints of $\Lambda_5$.
We define the triangle $\Delta_L$ for other
values of $L$ in the same way.
In [{\bf CS\/}, \S 3] we prove that
$\Lambda_L \subset \Delta_L$ and that
the interior of $\Lambda_L$ lies in the
interior of $\Delta_L$.   Finally,
we show that $\partial_0 N$ intersects
$\Delta_L$ only at the top vertex. 
These ingredients combine to show that
$\Lambda_L \cap \partial N=\emptyset$, and
this shows that $E(V)$ and
$E(W)$ cannot be equal.

\subsection{The Differential Equation}
\label{diffeq}

We will now go into more detail about
how the Bounding Triangle Theorem is proved.
Let $\ell=L/2$.
We consider the {\it backwards flow\/}
along the structure field $\Sigma$, namely
\begin{equation}
  \label{MAINDIFF}
  x'=-xz, \hskip 30 pt
  y'=+yz, \hskip 30 pt
  z'=x^2-y^2,
\end{equation}
with initial conditions $x(0)>y(0)>0$ and $z(0)=0$
chosen so that the point is in the loop level
set of period $L$. (We will often denote
these functions as $x_L$, etc.)
We let $g_t$ be the small symmetric flowline
whose endpoints are $(x(t),y(t),z(t))$ and
$(x(t),y(t),-z(t))$.  Then
$$\Lambda_L(t)=(a(t),b(t),0)=\Lambda_{g_t}.$$
Taking the derivative, we have
$$
(a',b',0)=\Lambda'_L(t)=\lim_{\epsilon \to 0} \frac{\Lambda_L(t+\epsilon)-\Lambda(t)}{\epsilon},
$$
$$
  \Lambda_L(t+\epsilon) \approx (\epsilon x, \epsilon y,\epsilon z) *
  (a,b,0) * ( \epsilon x, \epsilon y,-\epsilon z).
$$
The approximation is true up to order $\epsilon^2$ and
$(*)$ denotes multiplication in Sol.
A direct calculation gives
\begin{equation}
  a'=2x+az, \hskip 30 pt b'=2y-bz.
\end{equation}
The initial conditions are $a(0)=b(0)=0$.
(We will often denote these functions as $a_L$ and $b_L$.

\begin{lemma}
  \label{transform}
  For any $r \geq 0$ we have
  $$a(r)x(r)=\int_0^{r} 2x^2 dt, \hskip 30 pt
  b(r)y(r)=\int_0^{r} 2y^2 dt.$$
  \end{lemma}

\startproof
We have $(ax)'=2x^2$ and $(by)'=2y^2$.  Also
$a(0)=b(0)=0$.
Now we simply integrate.
\endproof

\begin{lemma}
\begin{equation}
\label{recip2}
\frac{b(0)}{a(0)}=\frac{b(\ell)}{a(\ell)}.
\end{equation}
\end{lemma}

\startproof
This comes from L'hopital's rule and the Reciprocity Lemma.
Let us take the opportunity to give a swift proof here.
(This is another proof of the Reciprocity Lemma in a special case.)
By two applications of Lemma \ref{transform}, we have
$$a(\ell)x(\ell)=\int_0^{\ell} 2x^2\ dt, \hskip 30 pt
b(\ell)y(\ell)=\int_0^{\ell} 2y^2\ dt.$$
But these two integrals are equal, by symmetry.
Hence $a(\ell)x(\ell)=b(\ell)y(\ell)$. Finally,
we have $x(0)=y(\ell)$ and
$y(0)=x(\ell)$ by symmetry.  Combining these
equations gives the result.
\endproof

The function $b(t)/a(t)$ has the same value
at $t=0$ and $t=\ell$.  To finish the proof, we just
have to show that  $b(t)/a(t)$ cannot have a
  local maximum.  This boils down to the fact that
  $ab''-ba''=2ab(y^2-x^2)$, a quantity which is negative
  for $t<\ell/2$ and positive for $t>\ell/2$.
These properties, together with the fact that $a'>0$, force
$\Lambda_L \subset \Delta_L$.  See [{\bf CS\/}, \S 3]
for more details.
\newline

\newpage

\section{More Information about the ODE}

In this chapter we further explore the
ODE we introduced in the previous chapter.
The results here
do not appear in [{\bf CS\/}].  However, they are
rather similar in spirit to some of the results there.
Lemma \ref{transform} above has a lot of juice in it, and
we want to squeeze some more out.  The estimates
here will be useful when we consider the
projections of the Sol spheres into the Euclidean plane
$\Pi_Z$.

\subsection{Bounding the Coordinates}

  In the next lemma, $2^*$ refers to a number which
  we can make as close as we like to $2$ by taking $L$ sufficiently large.
  This result says that the boundary of the yellow region in Figure 3.1
  asymptotes to the lines $X=2$ and $Y=2$.  

  \begin{lemma}
    \label{ABOUND}
    $b(\ell)<2^*$.
\end{lemma}

  \startproof
  From Lemma \ref{transform} and symmetry we have
  $$y(\ell)b(\ell)=\int_0^{\ell} y^2=\int_{0}^{\ell/2} (x^2+y^2) dt.$$
  The last equality follows from the fact that the
  function $t \to x^2(t)+y^2(t)$ is periodic with period $\ell/2$.
  Since $b(\ell) \sim 1$ for large $L$, it suffices to
  prove that the integral on the right approaches $1$ as
  $L \to \infty$.

We have
\begin{equation}
\int_0^{\ell/2} (x^2+y^2)dt=
\int_0^{\ell/2} (x^2-y^2)dt + 2\int_0^{\ell/2} y^2 dt.
\end{equation}
Now observe that
\begin{equation}
  \label{integrate00}
\int_0^{\ell/2}(x^2-y^2)dt =
\int_0^{\ell/2} z'\ dt = z(\ell/2)-z(0)=z(\ell/2) \sim 1.
\end{equation}

To finish the proof, it suffices to show that
\begin{equation}
  \label{integrate1111}
  \int_0^{\ell/2} y^2dt  \sim 0.
  \end{equation}
Let $\alpha$ be such that $(x(0),y(0),0)$ lies in the
same loop level set as $$(\alpha,\alpha,\sqrt{1-2\alpha^2}).$$
Then $y \leq \alpha$ on $[0,\ell/2]$ because $y$ is monotone
increasing on this interval.
Hence, by Equation \ref{periodXX} and some algebraic manipulation,
$$\int_0^{\ell/2} y^2 dt \leq L_{\alpha} \times \alpha^2 =
2\alpha  \sqrt{\frac{4\alpha^2}{1+2\alpha^2}} \times
{\cal K\/}\bigg(\frac{1-2\alpha^2}{1+2\alpha^2}\bigg).$$
 Setting $m=\frac{1-2\alpha^2}{1+2\alpha^2}$, we see that
\begin{equation}
\label{decay}
\int_0^{\ell/2} y^2dt \leq 2\alpha \sqrt{1-m} \times {\cal K\/}(m).
\end{equation}
As $L \to \infty$ we have $\alpha \to 0$ and $m \to 1$
and ${\cal K\/}(m) \sim -\log(1-m)/2$.  Hence, the
right hand side of Equation \ref{decay} tends to $0$
as $L \to \infty$.
\endproof

\noindent
{\bf Remark:\/}  We also have
 $a(\ell)b(\ell) \sim e^{\ell}$,
 as discussed in [{\bf CS\/}, \S 3.7] and also on
 [{\bf G\/}, p 75].

\begin{lemma}
\label{id1}
$a(L)=2b(\ell)$.
\end{lemma}

\startproof
We have 
$$a(L)x(L)=\int_0^{L} 2x^2 dt =
\int_0^{\ell} 2x^2 dt + \int_{\ell}^L 2x^2 dt=
2 \int_0^{\ell} 2y^2 dt = 2b(\ell)y(\ell).$$
Hence
$$a(L)x(L)=2b(\ell)y(\ell).$$
But $x(L)=y(\ell)$, so we can cancel these terms to
get the desired equality.
\endproof

\subsection{The Doubling Lemma}
\label{doubling}

It will be useful for us to consider flowlines which end on the plane $Z=0$.
These are the initial halves of symmetric flowlines.  Here
{\it doubling\/} refers to comparing the first half of a
symmetric flowline with the whole thing.

Let $\underline g_t$ denote the first half of
the flowline $g_t$; it connects the initial point of $g_t$ to the
midpoint of $g_t$.  Say that the coordinates of
$\Lambda_{\underline g_t}$ are $(\underline a(t),\underline b(t), \underline c(t))$.
The coordinate $\underline c(t)$ is typically nonzero, but we do not care about it.
Define
\begin{equation}
  \label{MAINDIFF2}
  \underline \Lambda_L(t)=(\underline a_L(t),\underline b_L(t)) \subset \R^2.
\end{equation}
In the next lemma we identity $\Pi_Z$ with $\R^2$.

\begin{lemma}[Doubling]
  $\underline \Lambda(t)= \frac{1}{2}\Lambda_L(t)$.
\end{lemma}

\startproof
We have
$$
(\underline a',\underline b',\underline c')=
\lim_{\epsilon \to 0} \frac{\underline \Lambda_L(t+\epsilon)-\underline \Lambda(t)}{\epsilon},
\hskip 20 pt
\underline \Lambda_L(t+\epsilon) \approx (\epsilon x, \epsilon y,\epsilon z) *
  (a,b,c).
  $$
  Taking the limit, we find that
  \begin{equation}
\underline a'=z + \underline ax, \hskip 30 pt
\underline b'=z - \underline bx.
\end{equation}
(Also $\underline c'=z$.
We have the same initial conditions $\underline a(0)=\underline b(0)$ as above.
Now notice that this solution to this equation is given by
$\underline a=a/2$ and $\underline b=b/2$.   
\endproof

\begin{corollary}
\label{id2}
$b(\ell)=\underline a(L)$.
\end{corollary}

\startproof
We combine the Doubling Lemma and Lemma \ref{id1} to get the equation
$b(\ell)=a(L)/2=\underline a(L)$.
\endproof

Now we give some applications, which show how the Doubling Lemma
and our asymptotics above give us some specific information about
some geodesic segments in Sol.
Let $f(r,L)$ denote the flowline of length $r$ on the loop level
set of period $L$ which ends at the point $(x,y,0)$ with $x>y$.
Let
$$\Upsilon_r(L)=\Lambda_{f(r,L)}.$$
This is the endpoint of the corresponding geodesic segment.
The {\it isochronal curve\/} $L \to \Upsilon_r(L)$, for $L \in [r,\infty)$  will be a central object
later in the paper.

 \begin{enumerate}
 \item We have $\Upsilon_r(r) \sim (2,e^{r/2}/2,0)$ by Lemma \ref{ABOUND} and the
   remark after Lemma \ref{ABOUND} and symmetry.

 \item We have $\Upsilon_r(2r) \sim (e^r/4,1,*)$ by Case 1, and symmetry, and
   the Doubling Lemma.

 \item We have $\Upsilon_r(4r) \sim (e^r/2,*,*)$.
 \end{enumerate}

We do not need Item 3 for any purpose, so we will be a bit sketchy with the proof.
   The small symmetric flowline of which $f(r,4r)$ is the first half starts
   at $(0,0,z)$ and ends at $(0,0,-z)$ for the appropriate choice of $z$.
   The corresponding geodesic segment $\gamma$ of length $2r$
   connects the origin to a point
   in $\Pi_Z$ and remains nearly tangent to $\Pi_Y$. Also, $\gamma$
   starts and ends nearly vertically. In fact, $\gamma$
   is asymptotic to the geodesic segment considered in Lemma \ref{fact1}.
   Thus, the first coordinate of the far endpoint of $\gamma$ is
   asymptotic to   $e^{r}$.   By the Doubling Lemma, the first coordinate of
   $\Upsilon_r(4r)$ is asymptotic to $e^r/2$.

   For what it is worth, the second
   coordinate tends to $0$ as $r \to \infty$.  To see this,
   note that the product of the first two coordinates of $\gamma$ is,
   by Lemma \ref{transform},
   $$\frac{4}{x(r)y(r)} \int_0^r x^2 dt \int_0^r y^2 dt \sim \frac{4}{x(r)y(r)}\int_0^r y^2 dt <
   \frac{4r y(r)^2}{x(r)y(r)}=4r.$$
   The asymptotic estimate $\sim$ is Equations \ref{integrate00} -- \ref{integrate1111}.
   The last inequality comes from the monotone increasing nature of $y(t)$ for $t \in [0,r]$.
   The equality comes from the fact that $x(r)=y(r)$.

\newpage

\section{The Euclidean Projection}

\subsection{Notation}

We introduce some notation that we use through out
the chapter.  We will consider some quantity $F$ that
depends on the variable $r$ or the variable $L$.
The statement $F<\zeta^*$ means that, for any
$\zeta^*>\zeta$ we can make $F<\zeta^*$ provided
that we take $r$ or $L$ sufficiently large. 

\subsection{The Cut Locus Image}

Figure 6.1 shows a more of the yellow region in Figure 3.1.
In this section we will prove a result which is
equivalent to the statement that the horizontal asymptote
of the boundary curve is the line $y=2$.
Similarly, the vertical asymptote is the line $x=2$.

\begin{center}
\resizebox{!}{2.5in}{\includegraphics{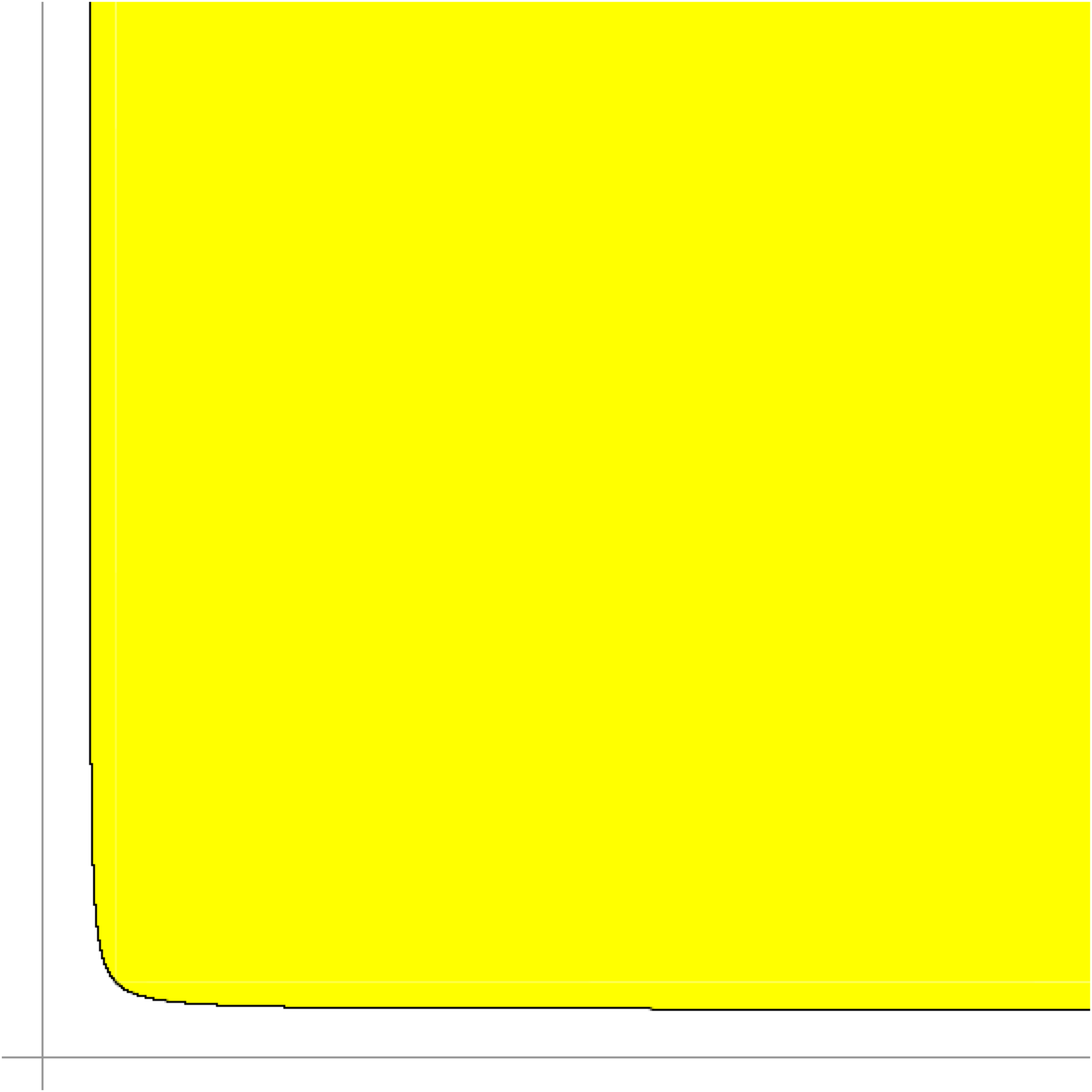}}
\newline
    {\bf Figure 6.1:\/} $\partial N$ in the positive quadrant.
\end{center}

Say that a vector is {\it positive\/} if all its
coordinates are positive. Say that a vector is
{\it distinguished\/} if its corresponding flowline is
contained in a small symmetric flowline with
the same initial point.  That is, the flowline
corresponding to the distinguished vector can be prolonged until it
is a small symmetric flowline.
In particular, distinguished vectors are small.

\begin{theorem}[Asymptotic]
  \label{asy}
If $V$ is a distinguished
positive vector whose corresponding flowline is
contained in a loop level set of period
$L$, then $E(V)=(a,b,0)$ has the property that
$b \in (0,2^*)$.
\end{theorem}

\startproof
First consider the special case when $V$ corresponds to
a small symmetric flowline.  Then
$b=b_L(t)$ for some $t \leq \ell$.
By the Bounding Triangle Theorem and Lemma \ref{ABOUND},
$b_L(t) \leq b_L(\ell)<2^*$.

Now consider the case when $V$ is an arbitrary
distinguished positive vector. There is some $\lambda \geq 1$
so that $\lambda V$ corresponds to a
small symmetric flowline.  The $X$ and $Y$ coordinates
of the curve $t \to E(tV)$ are increasing functions
because it is impossible for the geodesic associated to $V$
to be tangent to the hyperbolic foliations of Sol.
(Otherwise this geodesic would be trapped inside a
leaf of the foliation for all time.)  In particular,
the second coordinate of $E(V)$ is less or equal
to the second coordinate of $E(\lambda V)$, which
is in turn less than $2^*$ by the special case.
\endproof

\subsection{The Area Bound}

Our main goal is to show that
 $A_Z({\cal S\/}_r) <16^*e^r$.
Lemma \ref{conserve} below is
the main ingredient in the proof.
This result says that when $(a,b) \in \eta_Z({\cal S\/}_r)$
we have $\min(a^2b,ab^2)<2^* e^r$.
Figure 6.2 indicates the plausibility of this estimate.
Figure 6.2 shows the projection of (part of the positive sector of)
${\cal S\/}_5$ into
the plane $\Pi_Z$.  The small black arc of a hyperbola
is the projection of the singular set.
The outer black curve is
$\min(xy^2,x^2y)=e^5.$

\begin{center}
\resizebox{!}{2.3in}{\includegraphics{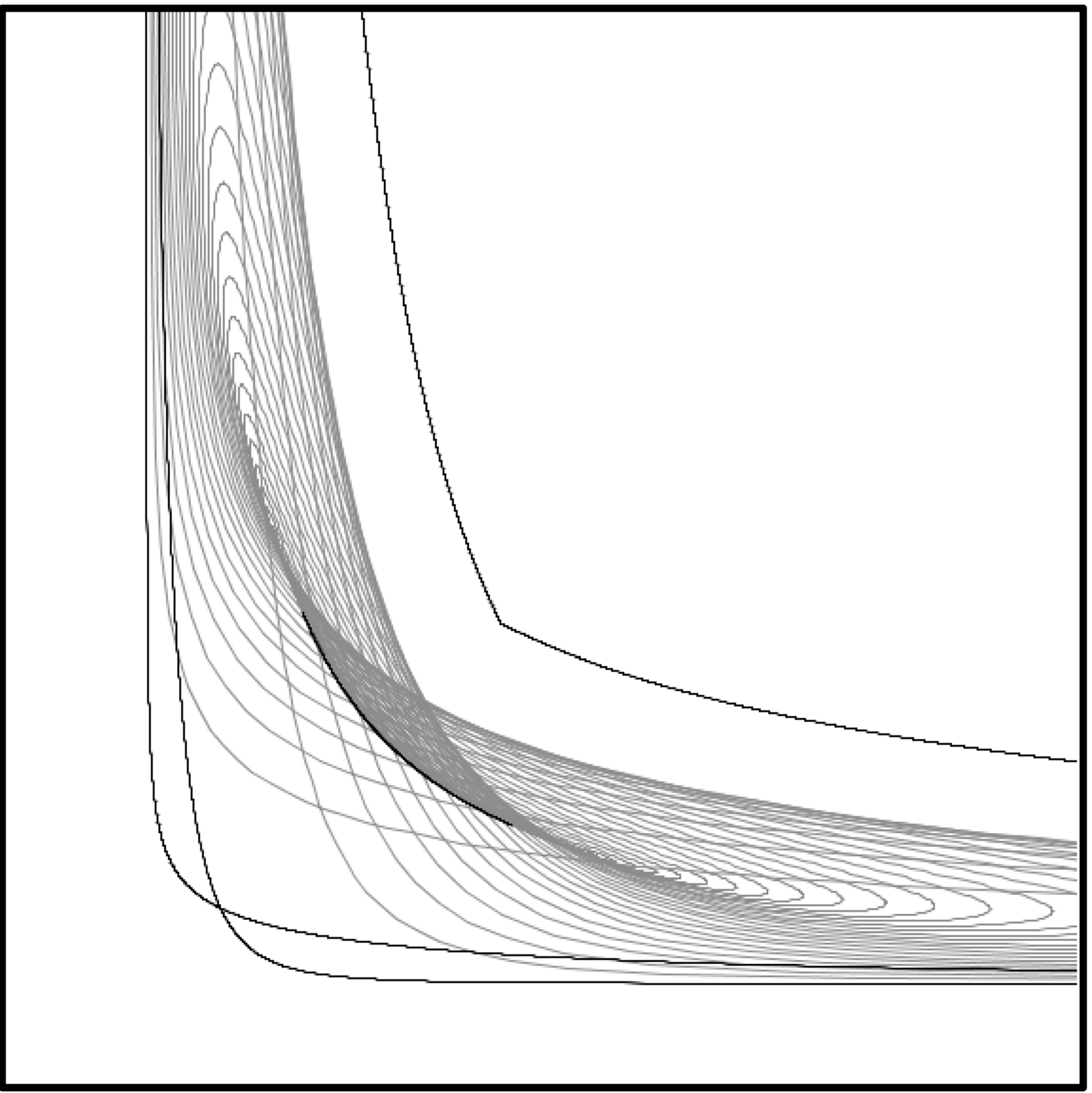}}
\newline
    {\bf Figure 6.2:\/} Projection into the plane $\eta_Z$.
\end{center}

\begin{lemma}
  \label{conserve}
Let $(a,b,c)=E(V)$, where $V$ is a small or perfect
vector of length $r$.
Then $\min(ab^2,a^2b)<2^* e^r$ and
$\max(a,b)<(1/2)^*e^r$.
\end{lemma}

\startproof
Note that as $r \to \infty$, the period of the loop level
set containing $V$ also tends to $\infty$.  This makes
the Asymptotic Lemma available to us.

Let $\gamma$ be the geodesic segment corresponding to $V$.
Let $g$ be the flowline corresponding to $V$. 
We can write $g=g_1|g_2$ where one of two things is true:
\begin{itemize}
\item $g_1$ is distinguished and $g_2$ is empty.
\item $g_1$ is small symmetric and $g_2$ is distinguished.
\end{itemize}
After interchanging the roles of $X$ and $Y$ if necessary,
we reduce to $2$ cases.
\newline
\newline
{\bf Case 1:\/}
Let $\gamma_1$ be the geodesic segment corresponding to
$g_1$, having endpoint $(a,b,c)$.
By the hyperbolic estimates, we know that
$\eta_Y(\gamma_1)$ lies in the hyperbolic disk $D_r$ in
$\Pi_Y$ centered at the origin.  By Lemma \ref{fact2},
we have $a<(1/2)^* e^r$.  We also have $b<2^*$by the Asymptotic Theorem.  Hence 
$ab^2<2^*e^r$.
We also see that $\max(a,b)<(1/2)^*e^r$.
\newline
\newline
{\bf Case 2:\/} 
Let $\gamma_1$ and $\gamma_2$ repectively be the
geodesic segments corresponding to $g_1$ and $g_2$.
Let $r_j$ be the length of $\gamma_j$.
Let $(a_j,b_j)$ be the projection to $\eta_Z$ of
the far endpoint of $\gamma_j$.

Let $\gamma_1'$ be the geodesic segment which
is the first half of $\gamma_1$, in terms of
length.  So, $\gamma_1'$ and $\gamma_1$ have
the same initial endpoint (the origin) but
$\gamma_1'$ has length $r_1/2$. 
Let $(a_1',b_1')$ be the far endpoint of
$\eta_Z(\gamma_1')$.  Once $r$ is large
enough we have the following:
\begin{itemize}
\item By Lemma \ref{fact1},  $a_1 \leq e^{r_1/2}-e^{-r_1/2}$.
\item By the Asymptotic Theorem,
$b_1<2^*$.
\item By Lemma \ref{fact2} and symmetry,
$b_2 \leq (e^{r_2}-e^{-r_2})/2$.
\item By symmetry and the Asymptotic Theorem,
$a_2<2^*$.
\end{itemize}
Combining these observations, we have
\begin{equation}
a \leq 2^*+ e^{r_1/2}-e^{-r_1/2}, \hskip 30 pt
b \leq 2^* + (e^{r_2}-e^{-r_2})/2.
\end{equation}
We have $r_1+r_2=r$.
We get right away that $\max(a,b)<(1/2)^*e^r$.

Now we consider $a^2b$.  Suppose first that
both $r_1$ and $r_2$ tend to $\infty$.
In this case $a<1^* e^{r_1/2}$ and
$b<(1/2)^*e^{r_2}$.  But then 
$a^2b<(1/2)^* e^r$.

When $r_1$ is bounded and
$r_2 \to \infty$,
$$a^2 b< (1/2)^* a^2 e^{r_2}=e^r \times 
\bigg(\frac{1}{2}+\frac{e^{-2r_1}}{2} - 2 e^{-3r_1/2} + e^{-r_1} 
+2 e^{-r_1/2}\bigg)<2^* e^r.$$
The last inequality follows from a bit of calculus.

When $r_2$ is bounded and $r_1 \to \infty$ we have
$a<1^* e^{(r-r_2)/2}$.  This gives
$$a^2 b< 1^* e^r \times \bigg(1-e^{-2r^2}+2e^{2r_2}\bigg)<2^* e^r.$$
This completes the proof.
\endproof

Once $r$ is sufficiently large, the sphere $S_r'$ consists entirely
of small and perfect vectors either contained in the
planes $\Pi_X$ and $\Pi_Y$ or else lying in loop level
sets whose period is so large that Lemma \ref{conserve} holds for them.
Lemma \ref{conserve} shows that the projection of the positive
sector of ${\cal S\/}_r$ lies in the region $\Omega_r$ defined by the
following inequalities.
\begin{equation}
  \label{slop}
  X,Y \in [0,(1/2)^*e^r], \hskip 30 pt
  \min(xy^2,yx^2)=2^*e^r.
\end{equation}
We set $x_0=y_0=(2^*e^r)^{1/3}$.
The region $\Omega_r$ is the union of the square
$[0,x_0] \times [0,y_0]$, whose area is $0^* e^r$,
and two other regions
which are swapped by reflection in the main diagonal
$x=y$.
One region lies underneath the graph $y=(2^*e^r/x)^{1/2}$
  starting at $x=x_0$ and ending at $x=(1/2)^*e^r$.  This region has area
\begin{equation}
\label{integrate1}
  \sqrt{2^*}  e^{r/2} \int_{x_0}^{(1/2)^*e^r} 
\frac{dx}{\sqrt x}<\sqrt{2^*} e^{r/2} \times 2 \times
  \sqrt{(1/2)^*e^r}<2^* e^r,
\end{equation}
  once $r$ is large.
  Hence $\Omega_r$ has area at most 
$4^* e^r$.  Recalling that $\Omega_r$
  contains the projection of the positive sector of ${\cal S\/}_r$, which
  is $1/4$ of the whole sphere, we see that
  $A_X({\cal S\/}_r)<16^* e^r$.
  \newline
  \newline
  {\bf Remark:\/}
The set $\eta_X({\cal S\/}_r)$ contains the region $G_r$  above the $X$-axis,
underneath the arc $\Upsilon_r[2r,4r]$ discussed in \S \ref{doubling},
and to the right of the line $x=e^2 \times e^{r/2}$.
Given our Embedding Theorem below, and the estimates in
\S \ref{doubling}, we see that the upper boundary of
$G_r$ is the graph of a decreasing function whose domain
has length $e^r/4^*$ and whose minimum is
asymptotic to $1$.  Thus $G_r$ contains a rectangle of area
$e^4/4^*$.  Hence $A_X({\cal S\/}_r)> (2/1^*) e^r$.

  \subsection{The Yin Yang Curve}
\label{yinyang}

Given $r$, we define the {\it yin-yang curve\/} $Y_r$ 
to be the set of points in $S_r'$ where the differential
$d(\eta_Z \circ E)$ is singular.  For $r \leq \pi \sqrt 2$ the
the curve $Y_r$ is connected. For $r>\pi \sqrt 2$, the
curve has $2$ disjoint components interchanged by the
map $(x,y,z) \to (y,x,-z)$.

\begin{center}
\resizebox{!}{2.4in}{\includegraphics{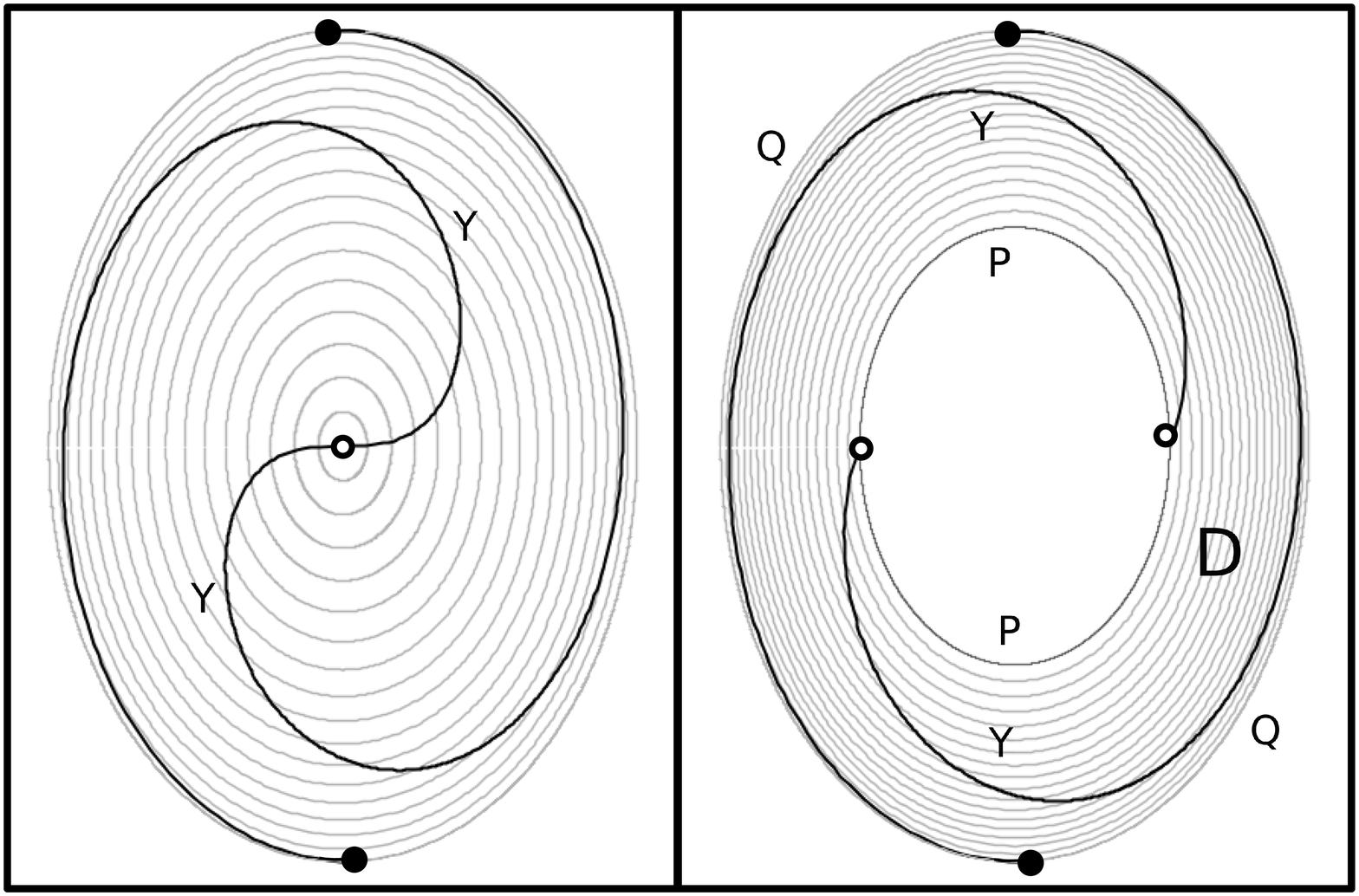}}
\newline
{\bf Figure 6.3:\/} The yinyang curves for $r=\pi \sqrt 2$ and $r=5$.
\end{center}

Figure 6.3 shows the yin-yang
curves for $r=\pi \sqrt 2$ and for $r=5$.  In Figure 6.3,
we are projecting the unit sphere in $\R^3$ onto the
plane through the origin perpendicular to the vector $(1,-1,0)$.
The loop level sets all project to ellipses having
aspect ratio $\sqrt 2$.  On the right side of 
Figure 6.3, the ellipse labeled $P$ is the set 
of perfect vectors on $S_5'$.  The ellipse
labeled $Q$ is the intersection of the
positive sector of the unit sphere with the
planes $X=0$ and $Y=0$.  Notice that
$P_r \cup Y_r \cup Q_r$ divides $S_r'$ into a union
of $2$ disks.  The map $\eta_Z \circ \Pi$ is
nonsingular on the interior of these disks and
hence a local diffeomorphism.  This is what
is important for our Projection Estimate 4.

Referring to Figure 6.3, the union
$Y_r \cup P_r \cup Q_r$ separates the
positive sector of $S_r'$ into $2$ components.
The map $(x,y,z) \to (x,y,-z)$ interchanges
these components. Let $D_r$ be either of these disks.
Below, we will describe more clearly which of
the two choices we take to be $D_r$.
Figure 6.4 shows $\eta_Z \circ E(D_5)$.  Essentially
this is ``half'' of Figure 6.2.

\begin{center}
\resizebox{!}{3.4in}{\includegraphics{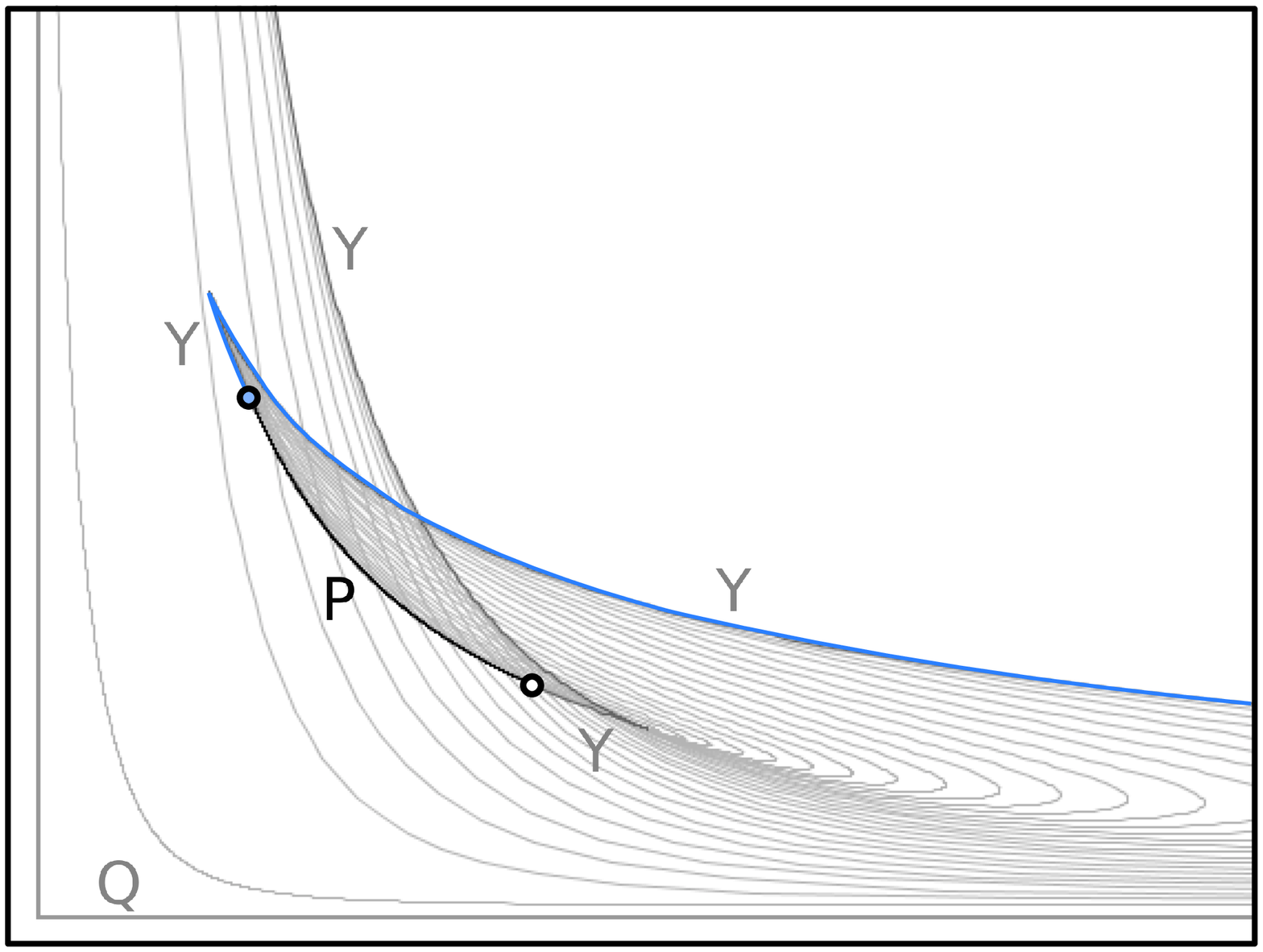}}
\newline
{\bf Figure 6.4:\/} The image $\eta_Z \circ E(D_5)$.
\end{center}

The region labeled $Y$ in Figure 6.3 is the image
$\eta_Z \circ E(Y_r)$.  We define
$\Upsilon_r=\eta_Z \circ E(Y_r^*)$, where $Y^*$ is
the component of $Y_r$ whose endpoint in
$\Pi_X$ is a point of the form $(x_r,y_r,0)$ with $x_r>y_r$.
We have drawn $\Upsilon_r$ in blue in Figure 6.4.

The map $L \to \Upsilon_r(L)$ is a smooth,
and indeed real analytic, map.  We say that
a {\it cusp\/} of this map is a point where
the map is not regular. So, away from the
cusps, $\Upsilon_r$ is a smooth regular curve.
Figure 6.3 suggests
that $\Upsilon_5$ just has a single cusp.
We prove the following result in the next chapter.

\begin{theorem}[Embedding]
  For $r$ sufficiently large, the curve $\Upsilon_r$
  has a single cusp, and negative slope away from a
single cusp. The cusp $\kappa_r=(a_r,b_r,c_r)$ satisfies the
bounds $a_r<2^*$ and $b_r<(e^2/2)^* e^{r/2}$.
  \end{theorem}

\subsection{The Multiplicity Bound}
\label{embedtheorem}

The image $\eta_Z \circ E(\partial D_r)$ is a piecewise
analytic loop.  We will show that this loop winds
at most twice around any point in the plane that it
does not contain.   
Referring to \S \ref{disk0}, we apply the
Disk Lemma to $h=\eta_Z \circ E$ and
$\Delta=D_r$.  This tells us that
$\eta_Z \circ E$ is at most $2$-to-$1$ on
$D_r$.  But then $\eta_Z \circ E$ is at most
$4$-to-$1$ on the positive sector of $S_r'$.
Hence $\eta_Z$ is at most $4$-to-$1$ on
the positive sector of ${\cal S\/}_r$.
Since the different sectors project into
$\eta_Z$ disjointly, 
we see that $\eta_Z$ is at most $4$-to-$1$ on
all of ${\cal S\/}_r$.  This establishes our
estimate $N_Z({\cal S\/}_r)=4$.

Now we turn to the analysis of the
image $\eta_Z \circ E(\partial D_r)$.
Define
\begin{equation}
  \Phi_r=E \circ \eta_Z(P_r).
\end{equation}
Also, let $R(x,y,z)=(y,x,-z)$.  The image
$\eta_Z \circ E(\partial D_r)$ is invariant
under $R$.  It is the union of $5$ analytic
arcs:
\begin{itemize}
\item An arc of the $X$-axis connecting the origin to the endpoint of $\Upsilon_r$.
\item $\Upsilon_r$.
\item $\Phi_r$.
\item $R(\Upsilon_r)$.
\item An arc of the $Y$-axis connecting an endpoint of $R(\Upsilon_r)$ to the origin.
\end{itemize}

Note that $\Upsilon_r \cup \Phi_r \cup R(\Upsilon_r)$
is a piecewise analytic arc that has its endpoints in
the coordinate axes and otherwise lies in the positive
quadrant.  The Embedding Theorem says that
$\Upsilon_r$ has negative slope, and is smooth and regular
away from a single cusp.
By symmetry, the Embedding Theorem also applies to
$R(\Upsilon_r)$.

The cusps serve as natural vertices for our loop,
so we make some new definitions which take
the cusps into account.
Let $\Phi_r^*$ denote the portion of
$\Upsilon_r \cup \Phi_r \cup R(\Upsilon_r)$
that lies between the two cusps.
Let $\Upsilon_r^*=\Upsilon_r - \Phi_r^*$.
The loop $\eta_Z \circ E(\partial D_r)$ has
the same $5$-part description as above,
with $\Upsilon_r^*$ and $\Phi_r^*$ used
in place of $\Upsilon_r$ and $\Phi_r$.

We will prove below that $\Upsilon_r^* \cup \Phi_r^*$ is
embedded.  By symmetry, $\Phi_r^* \cup R(\Upsilon_r^*)$ is
also embedded.  We will also prove that
$\Upsilon_r^*$ crosses the main diagonal -- the fixed
point set of $R$ -- exactly once.  This information
forces the schematic picture of
$\eta_Z \circ E(\partial D_r)$ shown in Figure 6.5.

\begin{center}
\resizebox{!}{2.2in}{\includegraphics{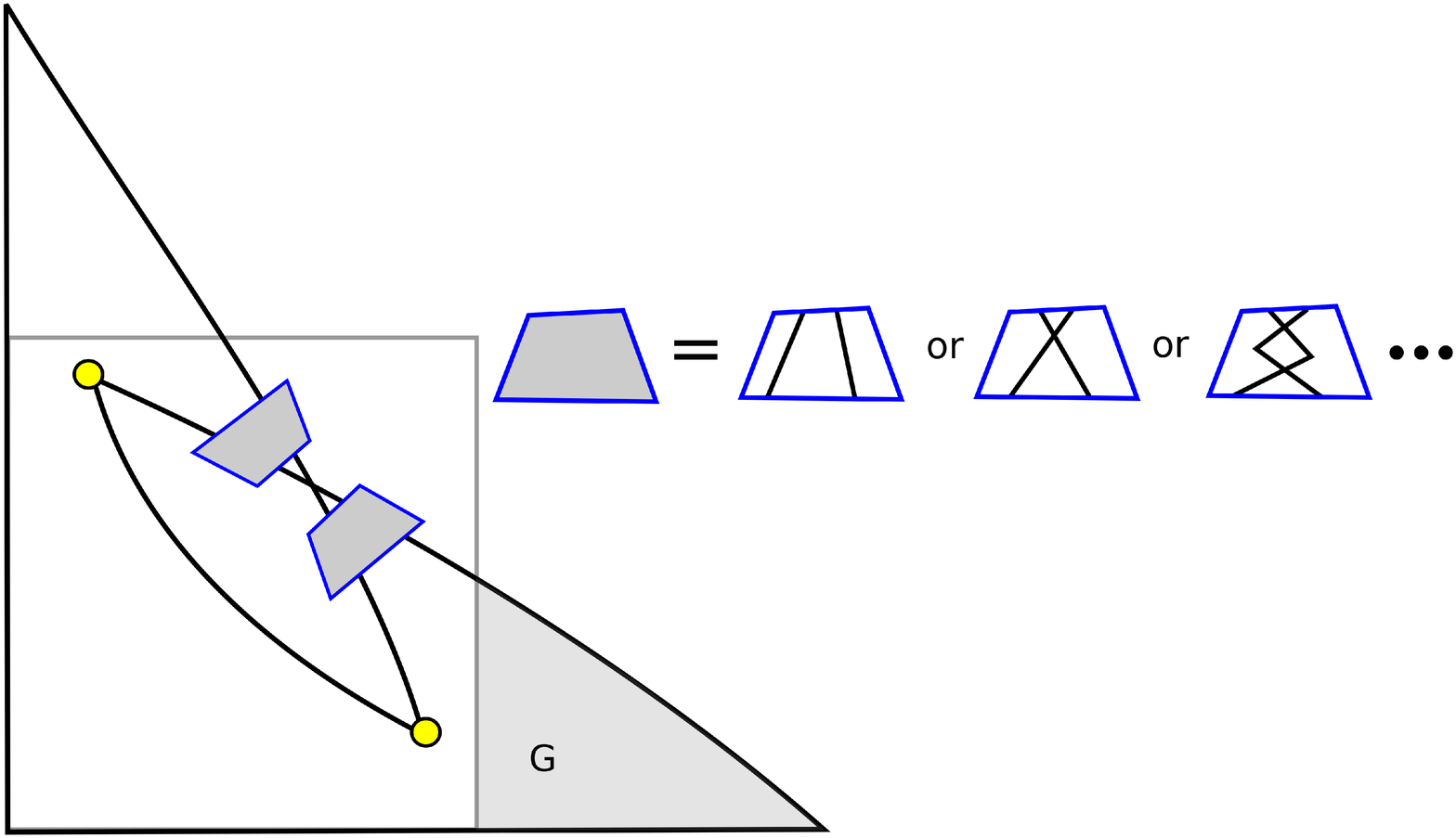}}
\newline
{\bf Figure 6.5:\/} Schematic picture of $\eta_Z \circ E(D_5)$.
\end{center}

Numerically, it seems that the first option occurs, and
that the two curves $\Upsilon_r^*$ and $R(\Upsilon_r^*)$
intersect exactly once.  We did not want to take the
trouble to establish this fact, given the already lengthy
nature of the paper.
In any case, the
information above establishes the fact that
  $\eta_Z \circ E(\partial D_r)$  winds
at most twice around any point in the plane
that does not lie in its image.
Applying the Disk Lemma
to the map $h=\eta_Z \circ E$ and the disk $D_r$ we see that
$h$ is at most $2$-to-$1$ on $D_r$.  But then $h$ is at most
$4$-to-$1$ on $D_r \cup I(R_r)=S_r^{++}$.  This completes the proof
of Projection Estimate 4,
modulo the properties of $\Upsilon_r^*$ and
$\Phi_r^*$.  We now turn to the task of
establishing the properties about the topology of
this planar loop.

Now we turn to the proof of Projection Estimate 5.
It follows from the negative slope
of $\Upsilon_r$ and $R(\Upsilon_r)$ that any
intersection between these two curves lies in
the square whose opposite corners are the two
cusps.   What we mean is that all the ``tangles''
shown in Figure 6.5 lie inside the lightly
shaded square.  Hence, $h$ is injective on the
portion of $D_r$ which maps outside this square.
Referring to the Embedding Theorem, the shaded
square is bounded by the lines $x=b_r$ and $y=b_r$,
where we know that $b_r<(e^2/2)e^{r/2}$.
The same analysis as done in connection with
Equation \ref{integrate1} shows that
$$A_{r,3}+A_{r,4}<K_1 e^{2r/3}+2I_r,$$ where
\begin{equation}
\label{integrate2}
I_r=  \sqrt{2^*}  e^{r/2} \int_{x_0}^{(e^2/2)^*e^{r/2}} 
\frac{dx}{\sqrt x}<K_2 e^{3r/4}.
\end{equation}
These bounds show that $A_{r,3}+A_{r,4}=0^*e^r$.
This is Projection Estimate 5.

\subsection{The Topology of the Boundary}
\label{bounds}

  \begin{lemma}
    For $r$ sufficiently large,
    the curve $\Upsilon_r^* \cup \Phi_r^*$ is embedded.
  \end{lemma}

  \startproof
Our argument refers to Figure 6.6.
  Let $\gamma$ be the portion of  $\Upsilon_r^* \cup \Phi_r^*$ 
  that lies above the (red) horizontal line $L_1$
  through the cusp of $R(\Upsilon_r^*)$.
  This point is the endpoint of $\Phi_L^*$.
  
\begin{center}
\resizebox{!}{2.5in}{\includegraphics{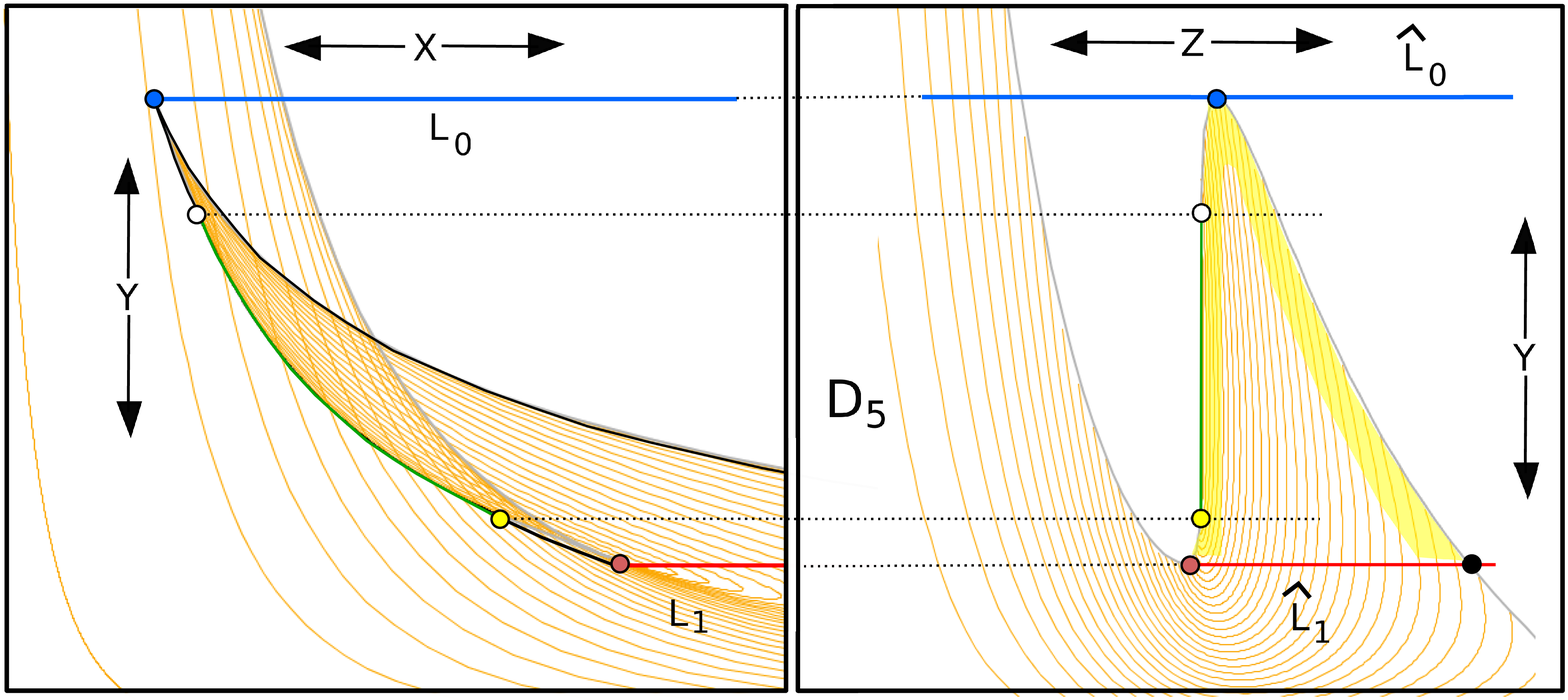}}
\newline
    {\bf Figure 6.6:\/} Projections of the relevant sets.
\end{center}

By the Embedding Theorem,
$\gamma$ has a single cusp, namely the cusp
of $\Upsilon_r$.  Let $\gamma_1$ and $\gamma_2$ be
the two arcs of $\gamma$ on either side of this
arc.  These two arcs have negative slope and
no cusps on them.   Hence all of $\gamma$ lies
between the horizontal line $L_0$ through the
cusp of $\Upsilon_r$ and the horizontal
line through the cusp of $R(\Upsilon_r)$.
These are the red and yellow horizontal
lines on the left side of Figure 6.6.

Since $\gamma_1$ and $\gamma_2$ have negative
slope and no cusps, they are each embedded.
We just have to see that $\gamma_1$
cannot intersect $\gamma_2$.  We will suppose
that there is an intersection and derive a
contradiction.

The portion of the Sol sphere ${\cal S\/}_r$
lying in the positive sector is the union
of two disks, $D_r$ and $R(D_r)$, where
$R$ is the isometry extending our reflection
in the main diagonal of $\Pi_Z$, namely
$I(x,y,z)=(y,x,-z)$.  The common boundary
of these disks contains a curve
$\widehat \gamma$ which projects to
$\gamma$.  We have
$\widehat \gamma=\widehat \gamma_1 \cup \widehat \gamma_2$.
On the right side of Figure 6.6, one of these
arcs connects the blue vertex to the black
vertex, and the other one connects the blue vertex
to the red vertex.  The right side of Figure 6.6
shows the projection into the $YZ$ plane.

Because no plane tangent to ${\cal S\/}_r$ at an
interior point of $\widehat \gamma$ is vertical, each
plane of the form $Y={\rm const\/}$. intersects
each of $\widehat \gamma_1$ and $\widehat \gamma_2$
exactly once.  In particular, this is
true to the planes $\widehat L_0$ and
$\widehat L_1$ which respectively project
to $L_0$ and $L_1$ on the left side of
Figure 6.6.

One of the two disks $D_r$ or $R(D_r)$
has the property that it lies locally
One of  the two disks -- say $D_r$ -- lies between
  $\widehat L_0$ and $\widehat L_1$ in a neighborhood
of $\widehat \gamma$.  We are taking about the
yellow highlighted region on the right side of
Figure 6.6.   The interior of $D_r$ is transverse to
  the plane $\widehat L_1$ because $\eta_Z$ is a
  local diffeomorphism on the interior of $D_r$.
  But this means that $D_r \cap \widehat L_1$
  contains a smooth arc $\widehat \beta$ which connects
  the endpoint of $\gamma_1$ to the endpoint
  of $\gamma_2$.  Let $\beta=\eta_Z(\widehat \beta)$.
  We note the following
  \begin{itemize}
  \item $\beta$ is contained in the line $L_1$.
  \item The endpoints of $\beta$ coincide.
  \item The interior of $\beta$ is a regular curve.
  \end{itemize}
  These properties are contradictory, because
  $\beta$ would have to turn around in $L_1$
  at an interior point, violating the regularity.
  This contradiction establishes the result that
  $\gamma$ is embedded.

  It remains to consider the portion of
  $\Upsilon_r \cup \Phi_r \cup R(\Upsilon_r^*)$
  that lies below the horizontal line $L_1$ through
  the cusp of $R(\Upsilon_r^*)$.  We label
  so that $\Phi_r \cup R(\Upsilon_r^*) \subset \gamma_2$.
  Given that $\gamma_2$ has negative slope we see
  that $\Phi_r \cup R(\Upsilon_r^*)$ lies entirely
  above $L_1$.  But this means that the portion
  of $\gamma_1$ below $L_1$ is disjoint from
  $\gamma_2$.  Finally, the portion of
  $\gamma_1$ below $L_1$ is disjoint from
  the portion of $\gamma_1$ above $L_1$ because
  $\gamma_1$ is regular and has negative slope.
  \endproof

  Define
  \begin{equation}
    \Upsilon_r^{**}=\Upsilon_r[r,2r]-\Upsilon_r^*.
  \end{equation}
  This is the subset of $\Upsilon_r[r,2r]$ that occurs after
  the cusp.  Given our result above, the only self-intersections
  on the curve $\eta_Z \circ D(\partial D_r)$ occur where
$\Upsilon_r^{**}$ and $R(\Upsilon_r^{**})$.  These are
analytic arcs of negative slope, and they are permuted
by the map $R$.  Hence, the can only intersect finitely
many times, and their intersection pattern must be as
in Figure 6.6.  This completes the proof of Projection Estimate 4,
and hence the Volume Entropy Theorem, modulo the proof of
the Embedding Theorem.

The rest of the paper is devoted to proving
the Embedding Theorem.

\newpage

\section{The Embedding Theorem}
\label{embedproof}

\subsection{The Isochronal Curves}

Let $\underline \Lambda_L=
(\underline a_L,\underline b_L)$, as in
\S \ref{diffeq}.  Let $E$ be the Riemannian exponential map
and let $\eta_Z$ be projection into the plane $\Pi_Z$.
We have $\Upsilon_r=\eta_Z \circ E(Y_r^*)$, where
$Y_r^*$ is the relevant component of $Y_r$.

\begin{lemma}
  $\Upsilon_r(L)=\underline \Lambda_L(r)$.
\end{lemma}

\startproof
Recall that $S_r'$ is the set of perfect vectors of length
$r$ contained in the positive sector of $\R^3$.
Define the {\it positive side\/} of $S_r' \cap \Pi_Z$
to be those vectors of the form $(x,y,0)$ with $x>y>0$.
By definition,
\begin{equation}
  \Upsilon_r(L)=\eta_Z \circ E(Y_r^*),
\end{equation}
where $Y_r^*$ is the component of the yin yang curve
$Y_r$ which ends in the positive side.  The
vectors $V \in Y_r$ are characterized by the property that
the differential $d(\eta_Z \circ E)$ is singular at points of $Y_r$.

The kernel of the projection map $\eta_Z$ is spanned
by the vector $(0,0,1)$.  So,
the differential $d(\eta_Z \circ E)$ is singular at
$V$ if and only if $dE$ maps the tangent plane to
$S_r'$ at $V$ to a plane which contains the vector $(1,0,0)$.
But then $dE(N_V)$ is orthogonal to $(0,0,1)$.
Here $N_V$ is normal to $S_r'$ at $V$.
But this means that the third coordinate of
$dE(N_V)$ is $0$.  Given the connection between
the Haniltonian flow on $S_r'$ and the geodesics,
this situation happens if and only if
the flowline associated to $V$ ends in the
plane $\Pi_Z$.

In short,
$Y_r$ consists of those small or perfect vectors
of length $r$ whose corresponding flowlines
end in $\Pi_Z$.  But these flowlines are then
the initial halves of symmetric flowlines
which wind at most twice around their
loop level sets.  The points in $Y_r^*$
are the initial halves of symmetric flowlines
whose midpoints lie on the positive side of $S_r'$.
Moreover, these symmetric flowlines wind at most twice
around their loop level sets and every
amount of winding, so to speak, from $0$ times to
$2$ times, is achieved.  So, by definition
$\Upsilon_r=\underline \Lambda(r)$.
\endproof

  We call $\Upsilon_r$ an
  {\it isochronal curve\/} because it
  computes all the solutions to the
  differential equation at the
  fixed time $r$.  Figure 7.1 shows
  part of $\Upsilon_5$.  The blue curves
  are the various curves $\underline \Lambda_L[0,L]$.
        
  \begin{center}
  \resizebox{!}{3.7in}{\includegraphics{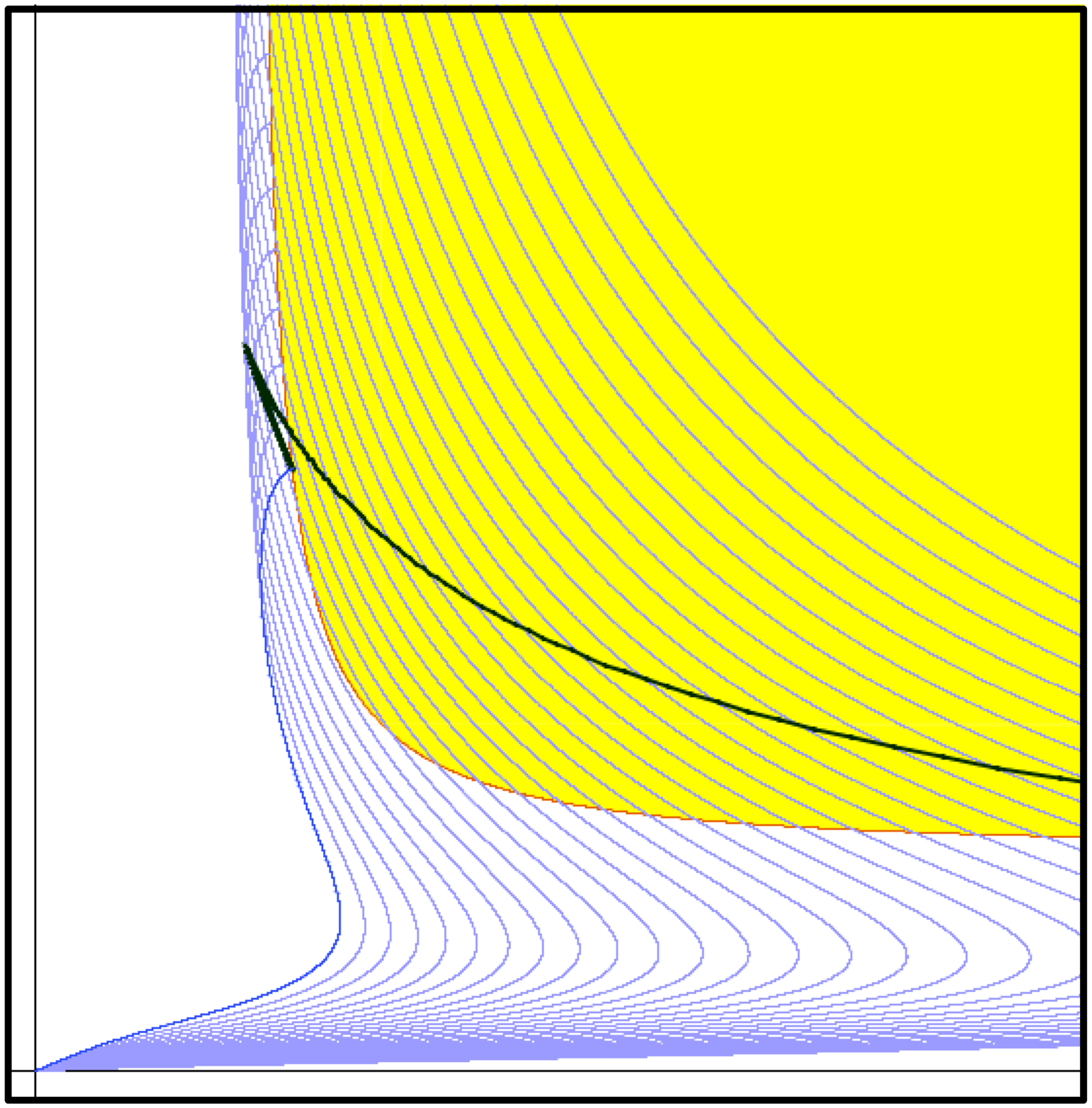}}
  \newline
  {\bf Figure 7.1:\/} The $\underline \Lambda_L$ curves and the
  initial part of $\Upsilon_5$.
  \end{center}

\subsection{The Tail End}
  
\begin{lemma}
  \label{tail}
  The curve $\Upsilon_r(2r,\infty)$ is smooth, regular,
    and embedded.
\end{lemma}

\startproof
The flowline corresponding to the point $\Upsilon_r(L)$  lies
on the loop level set of period $L$ and flows for time $r$.
If $L> 2r$ then the flowline travels less than halfway
around the loop level set.  Thus, the flowline is the
initial half of a small symmetric arc.
Let $S \subset \R^3$ denote the set of vectors corresponding
to small symmetric flowlines.
The map $E$ is a diffeomorphism on $S$, because $S$ consists
entirely of small vectors. This is part of the Cut Locus
Theorem from [{\bf CS\/}]. By the Doubling Lemma,
$$\Upsilon_r(2r,\infty)=\frac{1}{2} E(C_r),$$
where $C_r \subset S$ is a smooth regular curve, obtained
by dilating a suitably chosen arc of the yin yang curve by a factor of $2$.
Since $E$ is a diffeomorphism
on $S$, we see that $\frac{1}{2}E(C_r)$, is embedded.
\endproof

\subsection{The Slope}

Let $x_L$, etc. be the functions described
at the end of \S \ref{diffeq}.  These functions
satisfy the ODE
\begin{equation}
  x'=-xz, \hskip 15 pt
  y'=yz, \hskip 15 pt
  z'=x^2-y^2, \hskip 15 pt
  \underline a'=x+az, \hskip 15 pt
  \underline b'=y-bz,
\end{equation}
with initial conditions
\begin{equation}
  x_L(0)>y_L(0)>  z_L(0)=\underline a_L(0)=\underline b_L(0)=0,
  \hskip 15 pt (x_L(0),y_L(0),0) \in \Theta_L.
\end{equation}
Here $\Theta_L$ is the loop level set, in the positive
sector, having period $L$.
I am grateful to Matei Coiculescu for help with the
following derivation.

\begin{lemma}
\label{slopelemma}
  $\Upsilon_r$ has negative slope away from the cusps.
\end{lemma}

\startproof
This proof is a calculation with the ODE.
For any
relevant function $f$, the
notation $\dot f$ means
$\partial f/\partial L$.
We compute
$$(\underline{\dot a})'=
\frac{\partial }{\partial L}\frac{\partial \underline a}{\partial t} = \frac{\partial }{\partial L} (x+z\underline a)=
\dot x + \underline{\dot a} z + \dot z \underline a.$$
By the product rule
$$  \label{xa}
(x\underline{\dot a})' = x' \underline{\dot a} +
  x (\underline{\dot a})'=-zx \underline{\dot a} + x\dot z +
  x z \underline{\dot a} + \underline ax \dot z= x \dot x + \underline ax \dot z.
  $$
  This calculation, and a similar one, show that
\begin{equation}
  \label{yb}
  (x\underline{\dot a})' = x \dot x +\underline ax  \dot z, \hskip 30 pt
  (y\underline{\dot b})' = y \dot y -\underline by   \dot z
  \end{equation}
Since $x^2+y^2+z^2 \equiv 1$ we have
$x \dot x +  y \dot y + z \dot z  =0$.
Adding the Equations in Equation \ref{xa} and using
this relation, we find that.
$$(x \underline{\dot a} + y \underline{\dot b})'=
x \dot x + y \dot y + (ax-by)\dot z =
x \dot x + y \dot y + z\dot z = 0.$$
Hence $x \underline{\dot a} + y \underline{\dot b}$ is a constant function.
Since $\underline a(0)=\underline b(0)=0$ for all $L$ we have
$\underline{\dot a}(0)=0$ and $\underline{\dot b}(0)=0$.  So, the
constant in question is $0$.  Therefore
\begin{equation}
      \label{Matei3}
      x\underline{\dot a} + y \underline{\dot b}=0
\end{equation}
Because $\Upsilon_r(L)=\underline \Lambda_L(r)$,
  the velocity of
  $\Upsilon_r$ at $L$ is
  $$(\underline{\dot a}_L(r),\underline{\dot b}_L(r)).$$
  The slope of $\Upsilon_r$ at $L$ equals $-x_L(r)/y_L(r)$ by Equation
  \ref{Matei3}, provided that the velocity is nonzero.
  Since $x,y$ are everywhere positive,
  the slope of $\Upsilon_r$ is
  negative whenever the velocity is nonzero.
  \endproof

  \subsection{The End of Proof}

In the next chapter we prove the following result.
\begin{lemma}[Monotonicity]
  \label{estimate}
  If $L$ is sufficiently large then
  there is some $t_L \in (L-1,L)$ such that
  the function $\underline{\dot b}_L$
  is negative on $[L/2,t_L)$ and positive on $(t_L,L]$.
  Moreover, the function $L \to t_L$ is monotone increasing.
\end{lemma}

\begin{corollary}
  $\Upsilon_r[r,2r]$ has exactly one cusp.
\end{corollary}

\startproof
The curve $\Upsilon_r$ has a cusp at $L$ if and only if its velocity
$$(\underline{\dot a}_L(r),\underline{\dot b}_L(r))$$ vanishes.
  By Lemma \ref{Matei3}, one coordinate of the velocity vanishes
  if and only if the other one does.  So, $\Upsilon_r$ has a cusp at $L$
  if and only if $\underline{\dot b}_L(r)=0$.
    Suppose then that $\Upsilon_r[r,2r]$ has
  more than one cusp.  Then there are at
  least two pairs $(r,L_1)$ and $(r,L_2)$ such that
  $\underline{\dot b}_{L_1}(r)=0$ and  $\underline{\dot b}_{L_2}(r)=0$.
  This means that $r=t_{L_1}=t_{L_2}$.
  But this contradicts the Monotonicity Lemma.
  Hence $\Upsilon_r$ has at most one cusp.
  
  Since the map $L \to t_L$ is unbounded, each sufficiently
  large $r$ lies in its image.  But this means that for
  sufficiently large $r$ there is some $L \in (r,r+1)$ such that $r=t_L$.
  This means that $\Upsilon_r$ has a cusp at $L$.
  Hence, once $r$ is sufficiently large, $\Upsilon_r$ has
  exactly one cusp.
  \endproof
  
  This completes the proof of the Embedding Theorem, but
  there is one more remark we want to make.  The cusp of
  $\Upsilon_r$ occurs at some $L \in (r,r+1)$.  This explains
  why, in Figure 7.1, the cusp appears all the way to the
  left, near the end of $\Upsilon_5$.

  The next two chapters are devoted to the proof of the
  Monotonicity Lemma.

\newpage

\section{The Vanishing Point}

\subsection{Auxiliary Functions}

In this chapter we will prove the first half of the
Monotonicity Lemma.  That is, we will show that
once $L$ is sufficiently large
there is a point $t_L \in [L-1,L]$ such
$\underline {\dot b_L}$ vanishes at $t_0$,
and this is the only vanishing point. We sometimes set
$\dot f = \partial f/\partial L$ and
$f'=\partial f/\partial t$.
We introduce the following functions.

\begin{equation}
  \label{mainfunctions}
  Z = \dot z, \hskip 15 pt
  X=\dot x/x, \hskip 15 pt
  Y=\dot y/y, \hskip 15 pt
  B=\dot b/b, \hskip 15 pt
\end{equation}
Since $a,b,x,y>0$ these functions respectively  have
the same signs as $\dot z, \dot x, \dot y, \dot b$.

  \begin{center}
  \resizebox{!}{3.1in}{\includegraphics{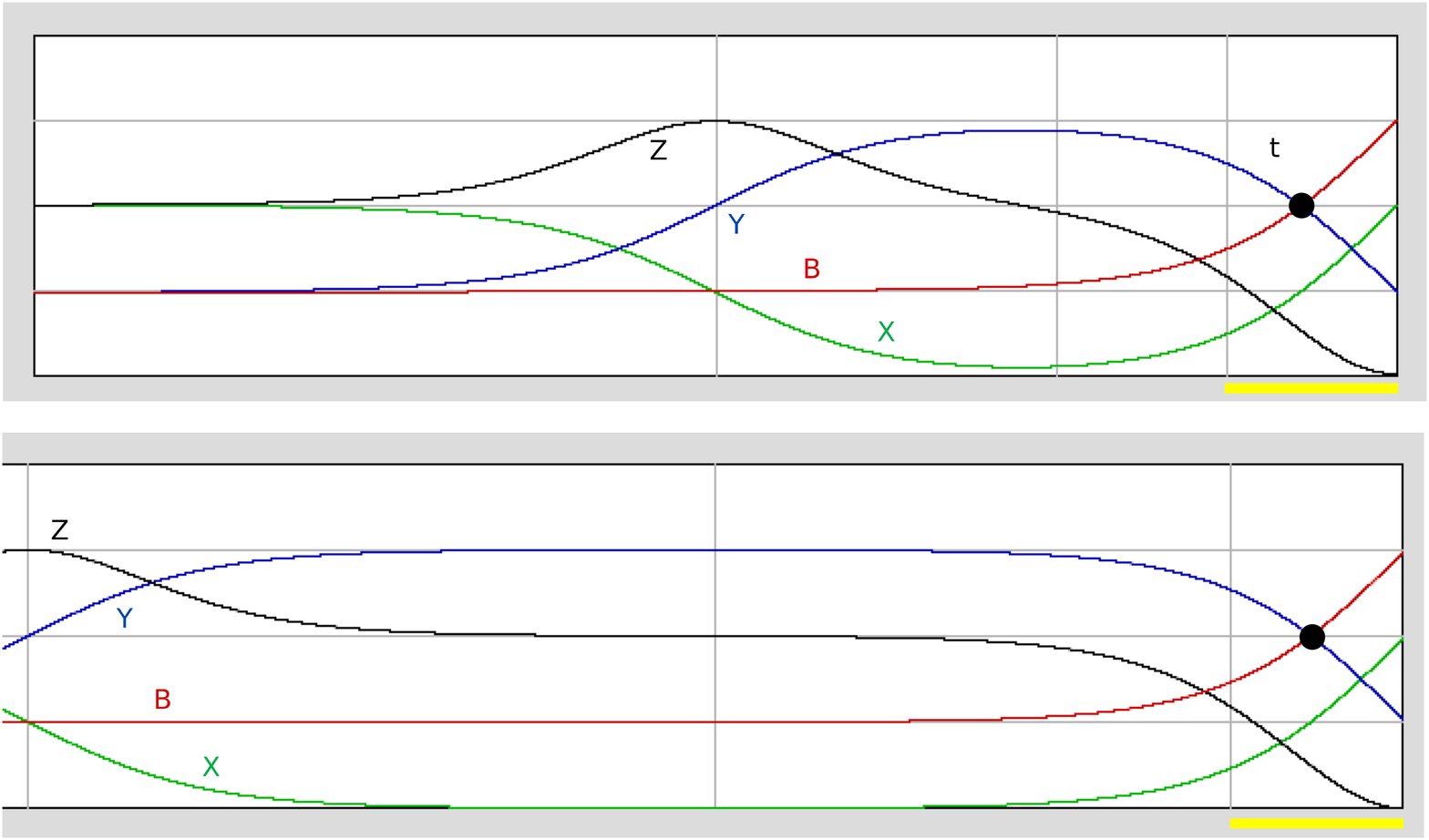}}
  \newline
      {\bf Figure 8.1:\/} Numerical plots
   \end{center}

The top half of
Figure 8.1 shows numerical plots of these functions at $L=8$.
The bounding box in the picture is $[0,8] \times [-1,1]$.
The bottom half shows the plot of the functions at the value $L=16$.
This time we are just showing the right half of the plot.  Notice that
in the interval $[L-1,L]$ the plots line up very nicely.
The black dot in both cases is the point $(t_L,0)$.  One half
of the Monotonicity Lemma establishes that $t_L$ is uniquely
defined.  The second half shows that $t_L$ increases monotonically.
The intuition behind the second half of the result is that
the pictures in $[L-1,L]$ stabilize, so that $t_L=L-s_L$
where $\partial s_L/\partial L$ is approximately $0$ for large $L$.

\subsection{Differentiation Formulas}

In this section we derive the following formulas.
\begin{enumerate}
\item  $X' =  -Z$
\item  $Y' =  +Z$
\item  $Z'=- 4y^2Y -2z Z.$
\item $B'= \frac{y}{\underline b}(Y - B) - Z.$
\item $Z''=(-2+6z^2)Z$.
\end{enumerate}

\noindent
(1) From the relation $x'=-xz$ we get
$(\dot x)'=-\dot x z - \dot z x$.
    Then we get:
$$X'=\bigg(\frac{\dot x}{x}\bigg)'=\frac{- \dot x z -  \dot z x}{x} - \frac{\dot x x'}{x^2}= - \dot z - \frac{\dot x z}{x}+\frac{\dot x z}{x}=-\dot z=-Z$$
(2) From the relation $y'=yz$ we get
$(\dot y)'=\dot y z + \dot z y$.
Then we get:
$$Y'=\bigg(\frac{\dot y}{y}\bigg)'=\frac{\dot y z +  \dot z y}{y} - \frac{\dot y y'}{y^2}=  \dot z + \frac{\dot y z}{y}-\frac{\dot y z}{y}=\dot z=Z$$
(3) From the relations $z'=x^2-y^2$ and $x^2 = 1-y^2-z^2$ we get
$$
Z'=2x\dot x + 2y \dot y=2x^2X-2y^2Y, \hskip 30 pt
2x^2 X = -2y^2Y -2zZ.$$
Substitute the second relation into the first to get the formula above. \newline \newline
(4) Note that $\underline{\dot b}/\underline b = \dot b/b$ because $\underline b=b/2$. So, we work
with $b$ for ease of notation.
From the relation $b'=2y - bz$ we get $(\dot b)'=2\dot y - \dot z b  -  \dot zb$. Then:
$$B'=\bigg(\frac{\dot b}{b}\bigg)'=
\frac{2 \dot y - \dot z b- \dot b z}{b}-\frac{\dot b b'}{b^2}=
\frac{2 \dot y - \dot z b- \dot b z}{b}-\frac{2\dot b y - \dot bbz}{b^2}=$$
$$\frac{2\dot y}{b} - \frac{2 \dot b y}{b^2} - \dot z =
\frac{2y}{b}(Y-B)-Z=\frac{y}{\underline b}(Y-B)-Z.
$$
(5)
We first work out that
\begin{equation}
  \label{diffZ}
  z''=-2z+2z^3.
\end{equation}
We then compute
$$Z'' = \frac{\partial z''}{\partial L}=
\frac{\partial(-2z+2z^3)}{\partial L} =
(-2+6z^2)Z.$$

\subsection{The Formula for B}
\label{Matei}

Here is a formula for $B$ in terms of the other quantities.
Matei Coiculescu found this formula and our derivation
follows his ideas.

\begin{equation}
  \label{nice}
B=\frac{x \underline a}{2y\underline b} X - \frac{1}{2} Y  - \frac{1}{2y\underline b} Z.
\end{equation}

We will work with
  $a$ and $b$ rather than $\underline a$ and $\underline b$ until the end.
  We have
$$
    (ax-by)'=2x^2-2y^2=2z', \hskip 20 pt a(0)=b(0)=0.
$$
  Integrating, we get
  \begin{equation}
    \label{Matei2}
    ax - by = 2z.
  \end{equation}

  Given that $a=2\underline a$ and $b=2\underline b$, Equation
  \ref{Matei3} from the previous chapter is equivalent to:
         \begin{equation}
      \label{Matei33}
      x\dot a + y \dot b=0.
       \end{equation}
  
Differentiating Equation \ref{Matei2} with respect to $L$, we have
\begin{equation}
  \label{Matei4}
x \dot a  + a \dot x -y \dot b  - b \dot y=2 \dot z.
\end{equation}  
Subtracting Equation \ref{Matei4} from Equation \ref{Matei33} we get
\begin{equation}
  \label{Matei5}
  2y \dot b - a \dot x  + b \dot y = -2 \dot z.
\end{equation}
Rearranging this, we get
$$
  \dot b = \frac{a\dot x - b \dot y - 2 \dot z}{2y}.
$$
When we make the substitutions
$$
  \underline a=a/2, \hskip 12 pt
  \underline b=b/2, \hskip 12 pt
  X=\dot x/x, \hskip 12 pt
  Y=\dot y/y, \hskip 12 pt
  Z=\dot z, \hskip 12 pt
$$
  we get Equation \ref{nice}.

  There is a similar formula for the
  function $A$, but we don't need it.

  \subsection{Elliptic Function Calculations}
  
Now we present the two calculations which we will use
in the proofs of
Lemmas \ref{Yasy} and \ref{Yasy2} below.
Define
\begin{equation}
  \label{uselater}
  L(y)=4f(y){\cal K\/} \circ g(y),
\end{equation}
where
\begin{equation}
  f(y)=\frac{1}{\sqrt{1+2y \sqrt{1-y^2}}} \hskip 15 pt
  g(y)=\frac{1-2y \sqrt{1-y^2}}{1+2y \sqrt{1-y^2}}.
\end{equation}
Next define
\begin{equation}
  Y_L=\frac{1}{y \times \frac{dL}{dy}}.
\end{equation}
We use the notation $f_L \sim g_L$ if $\lim_{L \to \infty} f_L/g_L=1$.

\begin{lemma}
  \label{uselater2}
  $Y_L \sim -1/2$ and $\frac{d}{dL} Y_L \sim 0$.
\end{lemma}

\startproof
We compute this in Mathematica, and we get a polynomial expression
terms of $\sqrt y$ and $K={\cal K\/} \circ g(y)$ and $E={\cal E\/} \circ g(y)$.  Taking the
series expansions of the coefficients, we find that
\begin{equation}
  \label{ASY1}
  Y_L = \frac{-1-5y+ \cdots}{\Delta}, \hskip 15 pt
  \Delta = (2+ 12 y+\cdots) E + (-4y - 24 y^2 + \cdots) K.
    \end{equation}
Given that $g(y) \sim 1-4y$, we see from Equations \ref{ASY} and
\ref{ASY1} that $Y_L \sim -1/2$.
Next we compute that
\begin{equation}
  \label{extradiff}
  \frac{\partial}{\partial y} Y_L = \frac{(48y +  \cdots)E + (-48 y + \cdots)K}{\Delta^2}.
\end{equation}
We see from Equations \ref{ASY0}, \ref{ASY} and \ref{extradiff} that  $\frac{\partial}{\partial y} Y_L \sim 0$.
\endproof

\subsection{Asymptotics}

If $f$ and $g$ are functions of $L$, we write
$f \sim g$ if $f/g \to 1$ as $L \to \infty$.
In this section we prove the following results.
$$(X_L(L/2),Y_L(L/2),Z_L(L/2),B_L(L/2)) \sim (-1/2,0,1/2,-1/2),$$
\begin{equation}
  \label{asy1}
(X_L(L),Y_L(L),Z_L(L),B_L(L)) \sim (0,-1/2,-1,1/2).
\end{equation}
These various features are already apparent in Figure 8.1.
Here we have written e.g. $X_L$ in place of $X$ to explicitly
indicate how the quantity depends on $L$.

\begin{lemma}
  \label{yequal}
  $Y_L(0)=Y_L(L)$.
\end{lemma}

\startproof
The limits we take for $Y_L(0)$ and $Y_L(L)$ respectively are
\begin{equation}
  Y_L(0)=\lim_{\epsilon \to 0} \frac{y_{L+\epsilon}(0)-y_L(0)}{\epsilon y_L(0)}, \hskip 30 pt
  Y_L(L)=\lim_{\epsilon \to 0} \frac{y_{L+\epsilon}(L)-y_L(L)}{\epsilon y_L(L)}.
\end{equation}
Note also that
$y_L(L)=y_L(0)$ and $y_{L+\epsilon}(L)=y_{L+\epsilon}(\epsilon)$.
Hence
\begin{equation}
  \label{newlim}
  Y_L(L)=\lim_{\epsilon \to 0} \frac{y_{L+\epsilon}(\epsilon)-y_L(0)}{\epsilon y_L(0)}.
  \end{equation}
But the map $t \to Y_{L+\epsilon}(t)$ has a local minimum at $t=0$ and
so $y_{L+\epsilon}(\epsilon)=y_{L+\epsilon}(0) + O(\epsilon^2)$.  Hence,
the limit in Equation \ref{newlim} equals $Y_L(0)$.
\endproof

\begin{lemma}
  \label{yequal2}
  $X_L(L/2)=Y_L(0)$.
\end{lemma}

\startproof
We have the relations
$$x_L(L/2)=y_L(0), \hskip 30 pt
x_{L+\epsilon}(L/2)=y_{L+\epsilon}(\epsilon/2)=y_{L+\epsilon}(0) + O(\epsilon^2).
$$
The last equality comes from the fact that the function
$t \to y_{L+\epsilon}(t)$ has its local minimum at $0$.
From these relations we have
$$\dot x_L(L/2)=\lim_{\epsilon \to 0} \frac{x_{L+\epsilon}(L/2)-x_L(L/2)}{\epsilon}= 
\lim_{\epsilon \to 0} \frac{y_{L+\epsilon}(0) + O(\epsilon^2)-y_L(0)}{\epsilon}=\dot y_L(0).$$
The lemma follows from this last equation and from $x_L(L/2)=y_L(0)$.
\endproof

\begin{lemma}
  \label{Yasy}
  \label{asymptote1}
  $Y_L(0) \sim -1/2$
\end{lemma}

\startproof
We set $y=y_L(0)$.
Using Equation \ref{periodXX} we get an explicit formula for $L$ in terms of $y$.
It is given by Equation \ref{uselater}.
By the Inverse Function Theorem we have
\begin{equation}
  Y_L(0)=\frac{1}{y \times \frac{\partial L}{\partial y}(y)}.
  \end{equation}
By Lemma \ref{uselater2}, we get
$Y_L(0) \sim -1/2$.
\endproof

\begin{lemma}
  \label{Xasy}
  $X_L(0) \sim 0$.
\end{lemma}

\startproof
We differentiate 
$x^2+y^2+z^2=1$ and use $z_L(0)=0$ to get
\begin{equation}
\label{sumsquare}
x_L^2(0)X_L(0) + y_L^2(0)Y_L(0)=0.
\end{equation}
We see that $X_L(0) \sim 0$ because
$x_L(0) \sim 1$ and
$y_L(0) \sim 0$ and $Y_L(0) \sim -1/2$ (a finite number).
\endproof

  \begin{lemma}
    \label{Zasy}
    $Z_L(L) \sim -1$ and $Z_L(L/2) \sim 1/2$.
    \end{lemma}

  \startproof
  We have
  $$z_L(L/2)=z_(L)=0, \hskip 20 pt
  z_{L+\epsilon}(L/2)=z_{L+\epsilon}(\epsilon/2), \hskip 20 pt
  z_{L+\epsilon}(L)=-z_{L+\epsilon}(\epsilon).$$
  Hence
$$
-Z_L(L) = \lim_{\epsilon \to 0} \frac{z_{L+\epsilon}(\epsilon)}{\epsilon}=
  \lim_{\epsilon \to 0} \frac{z_{L+\epsilon}(\epsilon)-z_L(\epsilon)}{\epsilon}
  +\lim_{\epsilon \to 0}\frac{z_{L}(\epsilon)}{\epsilon}.
$$
  The second limit on the right is just $z_L'(0)$.
  The first limit is $0$, because
  $$z_{L}(\epsilon) = z'_L(0) + O(\epsilon^2), \hskip 20 pt
  z_{L+\epsilon}(\epsilon)=z_{L+\epsilon}'(0) \epsilon + O(\epsilon)^2 =
  (z_L'(0)) \epsilon + O(\epsilon)^2.$$
  Hence $$Z_L(L)=-z_L'(0)=y_L^2(0)-x_L^2(0) \sim -1.$$
  The proof for $Z_L(L/2)$ is similar, and indeed
  $Z_L(L/2)=-(1/2)Z_L(L)$.
    \endproof

The rest of the relations for $X_L$ and $Y_L$
follow from the ones we have established above,
from Lemmas \ref{yequal} and \ref{yequal2}, and
the fact that $X_L+Y_L$ is a constant function.

\begin{lemma}
  \label{asymptote2}
  $B_L(L) \sim 1/2$ and $B_L(L/2) \sim -1/2$.
\end{lemma}

\startproof
Note that
$$(2\underline a_L(L/2),2 \underline b_L(L/2),
x_L(L/2),y_L(L/2))=(\underline b_L(L),\underline a_L(L),
y_L(L),x_L(L)).
$$
Hence, by the Reciprocity Lemma and
the Asymptotic Theorem,
$$2x_L(L/2) \underline a_L(L/2)=2y_L(L/2) \underline b_L(L/2) =
x_L(L) \underline a_L(L)=y_L(L) \underline b_L(L) \sim 2.$$
Equation \ref{nice} now gives us
$$B_L(L)=\frac{1}{2}X_L(L) - \frac{1}{2} Y_L(L)  - \frac{1}{2y_L(L)\underline b_L(L)} Z_L(L) \sim$$
\begin{equation}
  \label{nice3}
  \frac{1}{2} X_L(L) - \frac{1}{2} Y_L(L)  - \frac{1}{4} Z_L(L) \sim
  0  + (-1/2)(-1/2)  + (-1/4) (-1) =1/2.
\end{equation}
The proof for $B_L(L/2)$ is similar.
This completes the proof.
\endproof

\subsection{Variation of the Z Coordinate}
\label{ODE}

\begin{lemma}
  \label{Wron}
  Let $\pi \sqrt 2<L<M$. Then $z_M=z_L$ at most once on $(0,L)$.
\end{lemma}

\startproof
We already mentioned above that
$z''=-2z+2z^3$.
Let $z_L$ and $z_M$ be two solutions to this differential
equation. Consider the ratio $\phi=z_M/z_L$.
This quantity is positive on $(0,L/2]$.
By L'hopital's rule we can continuously
extend $\phi$ to $0$, and we have $\phi(0)>1$.

We have
\begin{equation}
\label{wron2}
\frac{d\phi}{dt} = \frac{W}{z_L^2}, \hskip 30 pt
W=z_L z_M' - z_M z_L', \hskip 30 pt
\end{equation}
If $z_L \not =0$ then $W$ and $\phi'$ have the same sign.
We compute
\begin{equation}
W'=z_L z_M'' - z_M z_L''=
2z_Lz_M(z_M^2-z_L^2).
\end{equation}

Suppose there is some smallest time $t \in (0,L/2]$ where $z_L(t)=z_M(t)$.
Note that $W' \not = 0$ on $(0,t)$.  Also, $W'>0$ on some small interval
$(0,\epsilon)$ because $z'_M(0)>z'_L(0)$.  Hence
$W'>0$ on $(0,t)$.  Hence $\phi$ is
increasing on $(0,t)$. In particular, $\phi(t)>1$.  Hence $z_M(t)>z_L(t)$, a contradiction.

If there is no time $t \in (L/2,L)$ where $z_M(t)=z_L(t)$,
then we are done.  Otherwise, let
$t_0$ be the smallest such time.
We have $z_L(t_0)=z_M(t_0)<0$.
Since $z_L(t_0-\epsilon)<z_M(t_0-\epsilon)$
for small $\epsilon>0$,
we have $z_M'(t_0) \leq z_L'(t_0)<0$.
Since these two functions are solutions of a
second order ODE, namely Equation \ref{diffZ},
they cannot have the
same initial conditions at $t_0$. Hence
$z_M'(t_0)<z_L'(t_0)<0$.

Let $\zeta_L=-z_L$ and $\zeta_M=-z_M$.  We
consider these functions on the interval
$(t_0,L)$.  These
are solutions of the same differential
equation, with initial conditions
$\zeta_L(t_0)=\zeta_M(t_0)$ and
$0<\zeta_L'(t_0)<\zeta_L'(t_0)$.
The same argument as above, applied to
$\zeta_L$ and $\zeta_M$, shows
$z_L(t)>z_M(t)$ for $t \in (t_0,L]$.
  \endproof

\begin{lemma}[Z variation]
$Z_L$ changes sign
  at most once on $[0,L]$ and $Z_L \geq 0$ on $[0,L/2]$.
\end{lemma}

\startproof
For the first statement, suppose
there are $3$ points $t_1,t_2,t_3$ where
$Z_L(t_i)$ for
$i=1,2,3$ alternates sign.  But then
for $\epsilon$ sufficiently small,
the difference
$z_{L+\epsilon}(t_i)-z_L(t_i)$ also alternates sign
for $i=1,2,3$. This contradicts Lemma \ref{Wron}.
So, there at most one $t_0 \in [0,L]$ where
$Z_L$ changes sign.  The second statment follows
from the analysis in Lemma \ref{Wron}, which
showed that $z_M(t)> z_L(t)$ when $L<M$ and $t \in (0,L/2)$.
\endproof

\subsection{Variation of the Y Coordinate}

\begin{lemma}[Y Variation]
Let $\delta_0=1/7$.
  \label{yvanish}
  If $L$ is sufficiently large then:
  \begin{enumerate}
\item $|Y_L|,|Z_L|<5$ on $[L-1,L]$.
  \item $Y_L(L)<-\delta_0$.
  \item $Y_L'<-\delta_0$ on $[L-1,L]$.
  \item $Y_L(L-1)>\delta_0$.
  \item $Y_L>0$ on $(L/2,L-1]$.
\end{enumerate}
\end{lemma}

\noindent
{\bf Proof of Statement 1:\/}
Let $\phi_L(t)=-Z_L(L-t)$.  This function satisfies the
differential equation
$$
  \phi_L(0)=-Z_L(L) \sim 1, \hskip 15 pt
  \phi'_L(0)=Z_L'(L) \sim 0, \hskip 15 pt
  \phi_L''(t)=(-2 + 6 z_L^2(t)) \phi_L(t).
$$
For the last equality we used the fact that $z_L^2(L-t)=z_L^2(t)$.
The solutions to this equation converge in the $C^{\infty}$ sense to
the solutions of the equation
\begin{equation}
  \phi(0)=1, \hskip 30 pt
  \phi'(0)=0, \hskip 30 pt
  \phi''=(-2 + 6 z^2) \phi.
  \end{equation}
Here $z$ satisfies $z''=-2z+2z^3$ with initial conditions $z(0)=0$ and $z'(0)=1$.
Since $\phi'' \in [-2,4] \phi$, and since $\cos(t \sqrt 2)>0$ on $(0,1]$, we have
\begin{equation}
\cos(t \sqrt 2) \leq \phi(t) \leq \cosh(2t).
\end{equation}
Since $\phi_L \to \phi$ we see that $|\phi_L|<2$ once $L$ is large.
Hence $|Z_L|<4$ on $[L-1,1]$.  Since $Y_L(L) \sim -1/2$ we have
$|Y_L|<5$ on $[L-1,L]$. \endproof
\newline
\newline
{\bf Proof of Statement 2:\/} This follows from the fact that $Y_L(L) \sim -1/2$.
    \newline \newline
{\bf Proof of Statement 3:\/}
Let $\phi_L$ and $\phi$ be the functions from the previous lemma.
For $t \in [0,1]$ we have (because $\phi$ is monotone decreasing)
$$
Y'_L(L-t)=Z_L(L-t)=-\phi_L(t)<\epsilon - \phi(t) \leq $$
\begin{equation} 
  \epsilon - \phi(1)  \leq
\epsilon - \cos(\sqrt 2) < -1/7
\end{equation}
once $\epsilon$ is sufficiently small.  We can arrange
this by taking $L$ sufficiently large. \endproof
\newline
\newline
{\bf Proof Statement 4:\/}
To estimate $Y_L(L-1)$ we note that
$$
Y_L(L-1) = Y_L(L) - \int_{L-1}^L Z_L\ dt =
Y_L(L) + \int_0^1 \phi_L\ dt > $$
\begin{equation}
- \epsilon -1/2 + \int_0^1 \phi\ dt=
-\epsilon  -1/2 + \frac{\sin(\sqrt 2)}{\sqrt 2}>1/7,
\end{equation}
once $\epsilon$ is sufficiently small. \endproof
\newline
\newline
{\bf Proof Statement 5:\/} Suppose this is false.
Since $Y_L(L/2) \geq 0$ and
$Y_L'(L/2)=Z_L(L/2)>0$ we see that
$Y_L$ is somewhere positive on $(L/2,L-1]$.
Also, $Y_L(L-1)>0$ and $Y_L(L)<0$.
If $Y_L=0$ somewhere else on $(L/2,L-1]$ then
  $Y_L$ switches signs at least $3$ times.
  But then $Z_L=Y_L'$ switches sign at least
  twice. This contradicts the Z Variation Corollary.
  \endproof

  \subsection{Variation of the B Coordinate}

\begin{lemma}[B variation]
  \label{close}
As $L \to \infty$, the function $Y_L+B_L$
converges to $0$ in the $C^1$ sense.
\end{lemma}

\startproof
We sometimes suppress the dependence on $L$ in our notation.
By Equation \ref{asy1} we have $|Y(L)+B(L)|<\epsilon/2$ if
$L$ is sufficiently large.  To finish the proof, it suffices to
show that $|B'+Y'|<\epsilon/2$ on $[L-1,1]$ for large $L$.
Combining our derivative formulas with the
preceding two results, we have
$$
|B'+Y'| \leq
\bigg|\frac{y}{\underline b}\bigg|\bigg(\max_{[L-1,L]} |Y|+|B|\bigg)
\leq 
\frac{10 y}{b}<10 y.$$
  For the last inequality, we note that $b(L/2)>1$ when $L$ is large and
  $b'>0$ on $[L/2,L]$.  Hence $b>1$ on $[L-1,1]$.
  As $L \to \infty$ the maximum value of $y$ on $[L-1,L]$ tends to $0$.
  \endproof

  \subsection{Uniqueness of the Vanishing Point}
\label{cusp00}

  In this section we prove the first half of the
  Monotoniticy Theorem.  That is, we show that
  $\underline{\dot b_L}$ vanishes
    exactly once on $(L/2,L)$.  This is the same
  saying that $B_L$ vanishes exactly once on
  $(L/2,L)$.  Our argument will establish the
stronger statement that $\underline{\dot b_L}$
vanishes somewhere for some $\lambda \in (L-1,L)$.
This establishes our assertion we made about
the cusp $\kappa_r$ just after the statement of the
Embedding Theorem in

\begin{lemma}
\label{cusplocation}
  $B_L$ vanishes exactly once in $[L-1,L]$, at an
  interior point $t_L$, and $B_L(L-1)<0$.
\end{lemma}

\startproof
Combining the Y Variation Lemma and the B Variation Lemma,
we get that $B_L(L-1)<0$ and $B_L(L)>0$ and $\dot B_L>0$ on $[L-1,L]$.
\endproof

    Let $t_1 \in [L/2,L-1]$ be the smallest value such that
    $Z \leq 0$ on $[t_1,L-1]$.
    
    \begin{lemma}
      $B_L<0$ on $[t_1,L-1]$.
    \end{lemma}

    \startproof
    Since $Z_L(L-1)<0$ we know that $t_1 \in [L/2,L-1)$.
    We introduce the function
    $\phi(t)=-B_L(L-1-t)$.
    We have $\phi(0)>0$ and
    $\phi'(t)=B_L'(L-1-t)$.
Hence
$$\frac{\partial \phi}{\partial t}=\frac{y_L(L-1-t)}{\underline b_L(L-1-t)} \times \bigg(Y_L(L-1-t) + \phi(t)\bigg) - 2 Z_L(L-1-t).$$
By the Y Variation Lemma and the definition of $t_1$ we have
$\phi'=\mu_1 \phi + \mu_2$
where $\mu_1$ and $\mu_2$ are non-negative functions on $[0,L-1-t_1]$.
Since $\phi(0)>0$ we have $\phi>0$ on $[0,L-1-t_1]$.
Hence $B_L<0$ on $[t_1,L-1]$.
\endproof

If $t_1=L/2$ this next lemma is vacuous.

\begin{lemma}
  $B_L<0$ on $(L/2,t_1)$.
  \end{lemma}

\startproof
Here is Equation \ref{nice} again:
$$
B=\frac{x \underline a}{2y\underline b} X - \frac{1}{2} Y  - \frac{1}{2y\underline b} Z.
$$
Our result here follows from Equation \ref{nice} and these inequalities on $[L/2,t_1)$:
\begin{itemize}
\item The quantities $x_L,y_L,\underline a_L,\underline b_L$ are all positive.
\item Since $X_L(L/2) \leq 0$ and $X_L'=-Z_L \leq 0$ on $(L/2,t_1)$, we have $X_L \leq 0$.
\item By the Y Variation Lemma, $Y_L>0$.
\item Since $Z_L$ changes sign only at $t_1$, and
  $Z_L(L)<0$, we have  $Z_L \geq 0$.
\end{itemize}
This completes the proof.
\endproof

\subsection{Bounds on the Cusp}

Now we prove the claim in the Monotonicity Theorem
concerning the location of the cusp
$$\kappa_r=(a_r,b_r,c_r).$$
We use the same notation established
at the beginning of \S 5.
We treat the bounds one at a time.

\begin{lemma}
$a_r<2^*$.
\end{lemma}

\startproof
Let $\Upsilon_r$ be the isochronal curve.
Let $\kappa_r'=(a_r',b_r',0)$ be the upper endpoint of the arc $P$.
We know from the Asymptotic Theorem and
from Lemma \ref{fact1} that $a_r'<2^*$.
The portion of $\Upsilon_r$ connecting $\kappa_r$ to $\kappa'_r$
has negative slope.  Hence $a_r<a_r'<2^*$. 
\endproof

\begin{lemma}
$b_r<(e^2/2)^* e^{r/2}$.
\end{lemma}

\startproof
By Lemma \ref{cusplocation}, we have
\begin{equation}\kappa_r \in \Upsilon_{\lambda}(L)
\end{equation}
for some
$\lambda \in (L-1,L)$.
Let $f_{r,\lambda}$ be the flowline corresponding to
$\kappa_r$.  From what we have just said,
the flowline $f_{r,\lambda}$ corresponding to $\kappa_r$ 
lies in the loop level set of period $\lambda$, winds
almost all the way around
its loop level set, and ends in the plane $Z=0$.

Consider the following perfect flowlines:
\begin{itemize}
\item $f_{\lambda}$ is the perfect symmetric flowline
which has the same ending point as $f_{r,\lambda}$.
Let $\Lambda_f \in \Pi_Z$ be the point corresponding
to $f_{\lambda}$.
\item $g_{\lambda}$ is the perfect symmetric flowline
which has the same initial point as $f_{r,\lambda}$.
Let $\Lambda_g \in \Pi_Z$ be the point corresponding
to $g_\lambda$.
\end{itemize}

By Lemma \ref{fact1} we have
\begin{equation}
\label{cusploc1}
\Lambda_f=(\alpha_f,\beta_f,0), \hskip 30 pt
\beta_f<(1/2)^* e^{\lambda/2}.
\end{equation}
Given the properties of concatenation, we have
some element $h \in {\rm Sol\/}$ whose third
coordinate lies in $(-1,1)$, such that
\begin{equation}
(\alpha_g,\beta_g,0)=\Lambda_g=h*\Lambda_f*h^{-1}.
\end{equation}
Conjugation by an element whose third coordinate
lies in $(-1,1)$ changes the first and second coordinates
by a factor of at most $e$.  Hence
\begin{equation}
\beta_g<(e/2)^* e^{\lambda/2}<(e^2/2)^* e^{r/2}.
\end{equation}
The
flowline $g_{\lambda}$ is the extension of
$f_{r,\lambda}$ by at most $1$ unit of flow.
Hence the distance from $\Lambda_g$ to $\kappa_r$ is less than $1$ unit.
Moreover, both points lie in the slab $|Z|<1$, where the
metric is boundedly close to Euclidean.  Hence, we get
 the same bound on $b_r$ as we got on 
$\beta_g$.  Hence $b_r<(e^2/2)^* e^{r/2}$.
The universally small
additive constant is just absorbed into the bound.
\endproof

\newpage

\section{Monotonicity of the Vanishing Point}

\subsection{The Proof Modulo Asymptotics}
\label{strategy}

We  know there is a unique $t_L \in (L-1,L)$
where $\underline{\dot b}_L$ vanishes.
In this chapter we show that $t_L$ is monotone increasing. Let
\begin{equation}
  \beta_L(t)=-B_L(L-t).
\end{equation}
There is a unique $s_L \in (0,1)$ such that
$\beta_L(s_L)=0$. In fact, $t_L=L-s_L$.
Since $\dot t_L=1-\dot s_L$, it suffices to show
$|\dot s_L|<1$ for large $L$.  We show $\dot s_L \sim 0$.
That is, $\lim_{L \to \infty} (\dot s_L)=0$.  Define
  \begin{equation}
  \delta(t)=\frac{y(L-t)}{\underline b(L-t)}, \hskip 30 pt
  \phi(t)=Y(L-t), \hskip 30 pt
  \zeta(t)=Z(L-t).
  \end{equation}
  Since $\beta'(t)=B'(L-t)$, we have
    \begin{equation}
      \beta'=\delta(\phi-\beta) - \zeta, \hskip 30 pt
    (\dot \beta)' = \dot \delta(\phi-\beta) +\delta(\dot \phi - \dot \beta) - \dot \zeta.
    \end{equation}
    The second equation comes from differentiating the first with respect to $L$.
    
    Let $\|f\|$ denote the sup of $f$ on $[0,1]$.
We will establish the following estimates below.
\begin{equation}
  \label{toprove}
  \|\delta\|, \|\dot \delta\|, \|\dot \phi\|, \|\dot \zeta\|, \dot \beta(0)  \sim 0.
\end{equation}
We also know from the Y Variation Lemma that
$\|\phi\|,\|\zeta\|<5$. Hence
\begin{equation}
  (\dot \beta)' =\epsilon_1 \dot \beta + \epsilon_2,
\end{equation}
where $\epsilon_1$ and $\epsilon_2$ are functions such that
$\|\epsilon_1\|, \|\epsilon_2\| \sim 0$.
Given our initial condition
$\dot \beta(0) \sim 0$, a standard comparison argument now says that
$\|\dot \beta \| \sim 0$.

By definition
\begin{equation}
  \beta_L(s_L)=0.
\end{equation}
Using implicit differentiation, we see that
\begin{equation}
  |\dot s_L| = \bigg|\frac{\dot \beta(s_L)}{\beta'(s_L)}\bigg|<8|\dot \beta(s_L)| \sim 0.
\end{equation}
The last inequality comes from the fact that $|\beta'|>1/8$ on $[L-1,L]$ once $L$ is large enough.
This proves the Monotonicity Lemma, modulo Equation \ref{toprove}.

\subsection{The Asymptotics}

We say that a function $f$ of $L$ is
{\it tame\/} if $df/dL \sim 0$.  If $f$ is analytic (i.e. not contrived in
an artificial way) and $f \sim {\rm const.\/}$
one might expect $f$ to be tame.  We verify that this is the case for a
number of quantities we have studied in the previous chapter.
We also point out, when relevant, how the given quantity relates
to the functions $\beta,\delta,\phi,\zeta$ introduced above.

\begin{lemma}
  \label{asymptote11}
  \label{Yasy2}
  $Y_L(L)$ is tame. Hence $\dot \phi(0) \sim 0$.
\end{lemma}

\startproof
Since $Y_L(0)=Y_L(L)$, it suffices to prove that $Y_L(0)$ is tame.
Lemma \ref{uselater2} tells us that
\begin{equation}
  \frac{d}{dy} Y_L(0) \sim 0,
\end{equation}
where $y=y_L(0)$.
But, since $Y_L(0)$ is asymptotic to a finite number, and $y_L(0) \sim 0$, we have
$\frac{dy}{dL} \sim 0$.
Therefore, {\it a fortiori} we have
$\frac{d}{dL}Y_L(0) \sim 0$.
\endproof

\begin{lemma}
$X_L(L)$ is tame.
\end{lemma}

\startproof
Since $X_L(0)=X_L(L)$, it suffices to prove that $X_L(0)$ is tame.
Differentiating Equation \ref{sumsquare} with respect to $L$, we get
\begin{equation}
  2 x_L^2(0) X_L^2(0) + x_L^2(0) \bigg(\frac{d}{dL} X_L(0)\bigg)+
  2y_L^2(0) Y_L^2(0) + y_L^2 \bigg(\frac{d}{dL} Y_L(0)\bigg)=0.
\end{equation}
Given that $x_L(0) \sim 1$ and $y_L(0) \sim 0$ and $Y_L(0) \sim -1/2$ and
that $Y_L(0)$ is tame, we see that $X_L(0)$ is also tame.
\endproof

  \begin{lemma}
    $Z_L(L)$ is tame.  Hence $\dot \zeta(0) \sim 0$.
    \end{lemma}

  \startproof
  We have
  $$\frac{d}{dL}Z_L(L)=\frac{d}{dL}(y_L^2(0)-x_L(0)^2)=
  y_L^2(0)Y_L(0) - x_L^2(0) X_L(0) \sim 0.$$
  The last equation comes from the fact that
  all quantities in the last expression are asymptotic to
  finite numbers and $y_L(0) \sim 0$ and $X_L(0) \sim 0$.
  \endproof

\begin{lemma}
  $Z'_L(L)$ is tame.  Hence $(\dot \zeta)'(0) \sim 0$.
\end{lemma}

\startproof
Since $z_L(L)=0$ we have 
\begin{equation}
\label{diffZprime}
Z'_L(L)=-4y_L^2(L) Y_L(L).
\end{equation}
Differentiating Equation \ref{diffZprime} with
respect to $L$ and using the product rule,
as above, we see that $\frac{d}{dL}Z'_L \sim 0$.
\endproof

\begin{lemma}
  $y_L(L) \underline b_L(L)$ is tame.
\end{lemma}

\startproof
We begin by proving an estimate that will come in at the end of
the proof. We claim that
\begin{equation}
  \label{YBOUND}
  \max_{[0,L/4]} |Y_L|<1.
  \end{equation}
  To see this, note that $Y_L'=Z_L$.  We also know that $Z_L \geq 0$ on $[0,L/2]$.
  Hence $Y_L$ is monotone increasing on $[0,L/2]$ and
  $Y_L(L/2) \sim 1/2$.  This establishes
  Equation \ref{YBOUND}.

By Lemma \ref{ABOUND},
we see that
\begin{equation}
  y_L(L)b_L(L)=2 z_L(L/4) + 2\phi_L, \hskip 20 pt
  \phi_L =\int_0^{L/4} y_L^2\ dt.
\end{equation}
We deal with these terms one at a time.
Referring to Equation \ref{periodXX} we have
$z_L(L/4)=\sqrt{1-2\alpha_L^2}$.
Here $\alpha_L \sim 0$.
This leads to  $\frac{d}{d\alpha}(Z_L(L/4)) \sim 0$.
A calculation like the one done in
Lemma \ref{asymptote1} shows that $|dL_{\alpha}/d\alpha| \sim \infty$.
Hence $d\alpha_L/dL \sim 0$.  But then, by the
chain rule, $\frac{d}{dL}z_L(L/4) \sim 0$.

We have
\begin{equation}
  \frac{d\phi_L}{dL}=\int_0^{L/4} \frac{\partial}{\partial L} (y_L)^2\ dt + \frac{1}{4}y_L(L/4).
\end{equation}
Referring to Equation \ref{periodXX} we have
$y_L(L/4)=\alpha$.  So, the same argument as for $Z_L(L/4)$ now shows that
$\frac{d}{dL} y_L(L/4) \sim 0$.
Finally, we have
$$\int_0^{L/4} \frac{\partial}{\partial L} (y_L)^2\ dt=
2\int_0^{L/4} y_L^2 Y_L\ dt \leq^* 2 \int_0^{L/4} y_L^2  \sim 0.$$
The starred inequality comes from Equation \ref{YBOUND}.
The final asymptotic result
comes from the proof of Lemma \ref{ABOUND}.
\endproof

All of the asymptotic results we have obtained so far feed into one final one.

\begin{lemma}
  $B_L(L)$ is tame.  Hence $\dot \beta(0) \sim 0$.
\end{lemma}

\startproof
As in Lemma \ref{asymptote2} have
$$B_L(L)=X_L(L) - \frac{1}{2} Y_L(L)  - \frac{1}{y_L(L)\underline b_L(L)} Z_L(L) \sim$$
\begin{equation}
  X_L(L) - \frac{1}{2} Y_L(L)  - \frac{1}{4} Z_L(L) \sim
  0  + (-1/2)(-1/2)  + (-1/4) (-1) =1/2.
  \end{equation}
We know that
$$X_L(L), Y_L(L), Z_L(L), y_L(L)\underline b_L(L)$$ are all tame.
Also, we know that $$y_L(L) \underline b_L(L) \sim 2,$$ by the
Asymptotic Theorem.
Using all this information, and the product and quotient rules for
differentiation, we see that  $B_L(L)$ is tame.
\endproof
  
  \begin{lemma}
    \label{ZZZ}
    $\|Z\| \sim 0$.
  \end{lemma}

  \startproof
  Our notation here is a bit funny. We mean
  to restrict our function $Z_L$ to the interval $[0,1]$
  and take its maximum.
    We have $Z(0)\sim 0$ and
    $|Z''| \leq 4|Z]$.  Since $Z>0$ on $(0,1]$,
    the same kind of comparison argument used in
    the proof of the Y Variation Lemma now shows that
    $\max_{[0,1]} |Z| \leq Z(0) \cosh(2) \sim 0$.
    \endproof

    Now we establish the remaining estimates from
    \S \ref{strategy}

    \begin{lemma}
      $\|\delta\| \sim 0$ and $\|\dot \delta\| \sim 0$.
    \end{lemma}

    \startproof
    The argument in the proof of the B Variation Lemma shows,
    incidentally, that $\|\delta\| \sim 0$.
        Using our derivative formulas, we have
    $$\dot \delta(t)= (Y(L-t)-B(L-t)) \delta(t).$$
  Combining the Y Variation Lemma and the B Variation Lemma we see that
  $|B|,|Y|<6$ on $[L-1,L]$ for large $L$.  Hence
  $\|\dot \delta\|<12 \|\delta\| \sim 0$.
    \endproof

  \begin{lemma}
    $\|\dot \zeta\| \sim 0$.
  \end{lemma}

  \startproof
  Define $\eta=\dot \zeta$.
 From the differential equation for $Z''$ and the fact that
 $z_{L-t}=z_t$ we get $\zeta''(t) = (-2 + 6z^2(t)) \zeta.$
 Differentiating with respect to $L$ and using the fact that
 the mixed partials commute, we see that
 $$\dot \zeta''(t) = (-2+6 z^2(t) \dot \zeta(t) + 12 Z(t)z^2(t)\zeta(t).$$
 The first term on the right lies in $[-4,4]\dot \zeta(t)$.
 The second term is at most
 $$12\|Z\| \max_{[L-1,L]}|Z|<60\|Z\| \sim 0.$$
 Here we have used Statement1 of the Y Variation Lemma
 and also the bound from Lemma \ref{ZZZ}.
 Putting these estimates together, we get
  \begin{equation}
   |\dot \zeta''| \leq 4 |\dot \zeta(t)| + \epsilon,
 \end{equation}
  where $\epsilon \sim 0$.
  We have already seen that
  $\dot \zeta(0) \sim 0$ and
  $(\dot \zeta)'(0) \sim 0$.
     The same kind of comparison argument as above
   now give us the desired bound on $\dot \zeta$.
 \endproof

  \begin{lemma}
    $\|\dot \phi\| \sim 0$.
  \end{lemma}

  \startproof
 We have $\phi'=-\zeta$.  Differentiating with
 respect to $L$ we get $$(\dot \phi)'=-\dot \zeta.$$
 We also have $\dot \phi(0) \sim 0$.
 We now integrate the 
 bound on $\|\dot \zeta\|$ to get the bound on
 $\|\dot \phi\|$.
  \endproof

 With these bounds, we complete the proof of the
 Monotonicity Lemma.
\newpage

\section{References}

\noindent
 [{\bf A\/}] V. I. Arnold, 
\textit{Sur la g\'eom\'etrie diff\'erentielle des groupes de Lie de dimension infinie et ses applications \`a l'hydrodynamique des fluides parfaits.} Ann. Inst. Fourier Grenoble, (1966).
\newline
\newline
[{\bf AK\/}] V. I. Arnold and B. Khesin, {\it Topological Methods in Hydrodynamics\/}, Applied Mathematical Sciences, Volume 125, Springer (1998) 
\newline
\newline
[{\bf B\/}] N. Brady, {\it Sol Geometry Groups are not Asynchronously Automatic\/}, Proceedings of the L.M.S., 2016 vol 83, issue 1 pp 93-119
\newline
\newline
[{\bf BB\/}] J. M. Borwein and P. B. Borwein, {\it Pi and the AGM\/},
Monographies et {\'E\/}tudes de la Soci{\'e}t{\'e} Math{\'e}matique du Canada, John Wiley and Sons, Toronto (1987)
\newline
\newline
[{\bf BS\/}] A.  B{\"o}lcskei and B. Szil{\'a}gyi, {\it Frenet Formulas and Geodesics in Sol Geometry\/},
Beitr{\"a}ge Algebra Geom. 48, no. 2, 411-421, (2007).
\newline
\newline
[{\bf BT\/}], A. V. Bolsinov and I. A. Taimanov, {\it Integrable geodesic flow with positive topological entropy\/},
Invent. Math. {\bf 140\/}, 639-650 (2000)
\newline
\newline
 [{\bf CMST\/}] R. Coulon, E. A. Matsumoto, H. Segerman, S. Trettel, {\it Noneuclidean virtual reality IV: Sol\/},
 math arXiv 2002.00513 (2020)
\newline
 \newline
[{\bf CS\/}] M. P. Coiculescu and R. E. Schwartz, {\it The Spheres of Sol\/}, submitted preprint, 2020
\newline
\newline
[{\bf EFW\/}] D. Fisher, A. Eskin, K. Whyte, {\it Coarse differentiation of quasi-isometries II: rigidity for Sol and Lamplighter groups\/},
Annals of Mathematics 176, no. 1 (2012) pp 221-260
\newline
\newline
[{\bf G\/}], M. Grayson, {\it Geometry and Growth in Three Dimensions\/}, Ph.D. Thesis,
Princeton University (1983).
\newline
\newline
[{\bf K\/}] S. Kim, {\it The ideal boundary of the Sol group\/}, J. Math Kyoto Univ 45-2
(2005) pp 257-263
\newline
\newline
[{\bf KN\/}] S. Kobayashi and K. Nomizu, {\it Foundations of Differential Geometry, Volume 2\/}, Wiley Classics Library, 1969.
\newline
\newline
[{\bf LM\/}] R. L{\'o}pez and M. I. Muntaenu, {\it Surfaces with constant curvature in Sol geometry\/},
Differential Geometry and its applications (2011)
\newline
\newline
[{\bf S\/}] R. E. Schwartz, {\it Java Program for Sol\/}, download (in 2019) from \newline
http://www.math.brown.edu/$\sim$res/Java/SOL.tar
\newline
\newline
[{\bf T\/}] M. Troyanov, {\it L'horizon de SOL\/}, Exposition. Math. 16, no. 5, 441-479, (1998).
\newline
\newline
[{\bf Th\/}] W. P. Thurston, {\it The Geometry and Topology of Three Manifolds\/}, \newline  Princeton University Notes (1978).
(See \newline http://library.msri.org/books/gt3m/PDF/Thurston-gt3m.pdf \newline
for an updated online version.)
\newline
\newline
[{\bf W\/}] S. Wolfram, {\it The Mathematica Book, 4th Edition\/},  Wolfram Media and Cambridge University Press (1999).

\end{document}